%% file: bozfer.tex
\documentclass[reqno, oneside]{amsart}
\usepackage{amsmath, amssymb, conformal}

\theoremstyle{plain}
\newtheorem{theorem}{Theorem}
\newtheorem{proposition}{Proposition}
\newtheorem{lemma}{Lemma}

\theoremstyle{definition}
\newtheorem{question}{Question}

\theoremstyle{remark}
\newtheorem*{remark}{Remark}

\def\fl(#1){\eqref{fl:#1}}
\def\thm#1{Theorem~\ref{thm:#1}}
\def\lem#1{Lemma~\ref{lem:#1}}

\def\prop#1{Proposition~\ref{prop:#1}}
\def\sec#1{\S{\rm \ref{sec:#1}}}

\def\op{\operatorname}
\def\Z{{\mathbb Z}}

\def\e{\varepsilon}
\def\f{\varphi}
\def\F{\op{F}}

\renewcommand\ge{\geqslant}
\renewcommand\le{\leqslant}
\renewcommand\~[1]{\widetilde{#1}}
\renewcommand\^[1]{\hat{#1}}
\renewcommand{\b}[1]{{\boldsymbol #1}}

\def\ol#1{\overline{#1}}
\def\inv{^{-1}}
\def\even{^{\bar 0}}
\def\odd{^{\bar 1}}
\def\<{\left<}
\def\>{\right>}
\def\({\left(}
\def\){\right)}
\let\ssm=\smallsetminus

\def\cff{\op{Coeff}}
\def\goth#1{{\mathfrak{#1}}}
\def\cal#1{{\mathcal{#1}}}

\def\:{\,\text{\Large \rm :}\,}

\newcommand{\ad}[2]{\b[\kern1pt #1,\, #2\,\b]}

\newcommand{\tr}[2]{\b\langle\kern1pt #1,\, #2\,\b\rangle}
\newcommand{\set}[2]{\left\{#1\ |\ #2\right\}}

\renewcommand\^[1]{\widehat{#1}}
\def\c{{\mathsf c}}

\begin{document}

\bibliographystyle{plain}

\author{Michael Roitman}
\thanks{Supported by Alfred P. Sloan Doctoral Dissertation
Fellowship. This research was partially conducted by the author for the Clay
Mathematics Institute.}

\address{M.S.R.I.,
1000 Centennial Dr.,
Berkeley, CA 94720.}  
\email{roitman@msri.org}

\title{Conformal subalgebras of lattice vertex algebras}

\date\today

\begin{abstract}
\input{abstract.tex}

\end{abstract}


\maketitle

\section*{Introduction}
\label{sec:intro}
\input{intro}

\section{Conformal and vertex algebras}
\label{sec:conformal}
\input{conformal.tex}
\input{examples.tex}

\input{vertex.tex}

\input{lattice.tex}

\section{Jordan triple systems and Tits-Kantor-Koeher construction}
\label{sec:tkk}
\input{tkk.tex}

\section{Bozon-Fermion correspondence}
\label{sec:bfc}
\input{fermion.tex}

\input{v0.tex}

\section{Conformal subalgebras of lattice vertex algebras}
\label{sec:solution}
\input{rootsys.tex}

\input{rank1.tex}

\input{rank2.tex}

\input{finite.tex}

\input{ears.tex}

\bibliography{conformal}

\end{document}

%% file: abstract.tex
In this paper we classify, under certain restrictions, 
all homogeneous conformal subalgebras $\goth L$ of a lattice 
vertex superalgebra $V_\Lambda$ corresponding to an integer lattice
$\Lambda$. We require that $\goth L$ is graded by an almost finite
root system $\Delta\subset \Lambda$ and that $\goth L$ is stable under
the action of the Heisenberg conformal algebra $\goth H\subset
V_\Lambda$.  We also describe the root systems of these
subalgebras. The key ingredient of this classification is an infinite
type conformal algebra $\goth K$ obtained by the Tits-Kantor-Koeher
construction from a certain  Jordan conformal triple system $\goth J$. 
We realize a central extension $\^{\goth K}$ of $\goth K$ inside the
fermionic vertex superalgebra $V_\Z$, thus extending the
bozon-fermion correspondence.

%% file: intro.tex
One of the first origins of vertex algebras was the explicit
constructions of representations of certain Lie algebras by means of
so-called vertex operators.  The first construction of this kind was
done by Lepowsky and Wilson \cite{lepwil}, who constructed a vertex
operator representation of the affine algebra $A_1^{(1)}$. Their work
was later generalized in \cite{kklw}.  Frenkel and Kac \cite{fk} and,
independently, Segal \cite{segal} constructed the basic
representations of the simply-laced affine Lie algebras using the 
so-called {\it untwisted vertex operators} as opposed to {\it twisted
vertex operators} of Lepowsky and Wilson.  Later vertex operators were
used to construct a large family of modules for different types of Lie
algebras, including all of the affine Kac-Moody algebras, toroidal
algebras and some other extended affine Lie algebras, see
e.g. \cite{billig,flm,gnos,kac1,mry,yamada} and references
therein. The advantage of vertex operator constructions is that they
are very explicit. They have yielded a lot of interesting results for
combinatorial identities, modular forms, soliton theory, etc. 

It seems to be a natural problem to describe all Lie algebras, or at
least a large family of Lie algebras, for which the vertex operator
constructions of representations work. Our first observation is that
in some of the cases described above the Lie algebras, whose
representation are constructed by vertex operators, correspond to 
{\it conformal algebras}, introduced by Kac \cite{kac2,kac_fd}, see
also \cite{primc,freecv}.  On the other hand, vertex operators give
rise to another algebraic structure, called {\it vertex algebras},
studied extensively in e.g. \cite{bor,fhl,flm,kac2}.  A vertex
operator construction of representations of Lie algebras amounts
sometimes to an embedding of a conformal algebra into a vertex
algebra generated by vertex operators, so that the vertex algebra
becomes an {\it enveloping vertex algebra} of these conformal
algebras.

In the present work we make the first step in describing the Lie
algebras representable by vertex operators. We classify the Lie
algebras that can be realized by the 
untwisted vertex operators of Frenkel-Kac-Segal.  The vertex algebra
generated by these vertex operators is also called {\it lattice vertex
algebra}, because its construction depends on a choice of an integer
lattice. In fact there is a functor that for every lattice
$\Lambda$ with an integer-valued bilinear form $(\,\cdot\,|\,\cdot\,)$  
gives a vertex superalgebra 
$V_\Lambda = \bigoplus_{\lambda\in\Lambda} V_\lambda$ graded by the 
lattice $\Lambda$, see \cite{dong,dl,flm,kac2}
and also \sec{lattice} of this paper. If $\Lambda$ is a 
simply-laced root lattice of a finite-dimensional simple Lie algebra,
then $V_\Lambda$ is a module over the corresponding affine Kac-Moody
algebra.  Lattice vertex algebras play an important role in different
areas of mathematics and physics, in particular the celebrated
Moonshine vertex algebra $V^\natural$, such that $\op{Aut}V^\natural$
is the Monster simple group, is closely related to the lattice vertex
algebra of certain even unimodular lattice of rank 24, called the
Leech lattice \cite{bor,flm}.

So the problem in our case takes the following form: describe all
conformal subalgebras  of the vertex operator superalgebra $V_\Lambda$
corresponding to an integer lattice $\Lambda$. 

Of course, we are not
interested in just any subalgebras of $V_\Lambda$, which are numerous
beyond any control, but rather in those subalgebras $\goth L \subset
V_\Lambda$ which can be
explicitly realized by vertex operators. A careful look at the results
cited in the first paragraph reveals that the requirements for 
$\goth L$ are the following:
\begin{itemize}
\item[\hskip15pt $\bullet$]
$\goth L = \bigoplus_{\lambda\in\Lambda} \goth L_\lambda$ is
homogeneous with respect to the grading by the lattice $\Lambda$. Here 
$\goth L_\lambda = \goth L\cap V_\lambda$. Let 
$\Delta = \set{\lambda\in\Lambda}{\goth L_\lambda\neq0}$ be the {\it root
system} of $\goth L$.
\item[\hskip15pt $\bullet$]
$|\Delta|<\infty$. In fact we can somewhat relax this requirement in
the case when $\Lambda$ is semi-positive definite, see \sec{rootsys}.
\item[\hskip15pt $\bullet$]
$\goth L$ is stable under the action of the Heisenberg conformal
algebra $\goth H\subset V_0$, see \sec{lattice}.
\end{itemize}
Yet another important thing is to study all possible root systems which
could be obtained in this way. 

Some unexpected results occur already in the case when the rank of
$\Lambda$ is 1, see \sec{rank1}. 
It turns out that besides well known realizations of
the Clifford and affine $\^{sl}_2$ algebras, one also gets the realization
of the $N=2$ simple conformal superalgebra, see e.g. \cite{kac98},  
and also of the conformal algebra $\goth K$ (or rather its central
extension $\^{\goth K}$),
obtained by the Tits-Kantor-Koeher construction from certain Jordan
conformal triple system $\goth J$, see \sec{tkk}. This exhausts all
possibilities for rank 1 lattices. In \sec{rank2} we extend this
result for the lattices of rank 2 and then in \sec{posdef}-\sec{ears}
we generalize the classification for the case of an arbitrary lattice.

We also classify all finite root systems of conformal
subalgebras of $V_\Lambda$, see \sec{finite}. Most of them occur if 
$\Lambda$ is positive definite. In this case all such indecomposable 
root systems $\Delta$ are in fact classical Cartan systems of type other than
$F_4$ and $G_2$. However the corresponding Lie algebras are affine
Kac-Moody only when $\Delta$ is simply-laced,  and then we are in the
situation of \cite{fk}.   
We also explain what goes on if the lattice $\Lambda$ is
semi-positive definite and outline the relation to the theory of extended
affine root systems (EARS), see \cite{aabgp}.

Another important motivation of the present work is the combinatorial
approach to vertex algebras. Once we know how to construct free vertex
algebras, see \cite{freecv}, we can consider presentations of vertex
algebras in terms of generators and relations. If a vertex algebra
$\goth V$ turns out to be an enveloping vertex algebra of a small conformal
algebra $\goth L$, then one can hope for getting a nice presentation 
of $\goth V$. For example, the Frenkel-Kac-Segal lattice vertex algebra of
an affine simply-laced Kac-Moody algebra is finitely presented. It
seems to be a very interesting problem to study presentations of other
lattice vertex algebras and also of the Frenkel-Lepowsky-Meurman Moonshine
vertex algebra $V^\natural$.

The paper is organized as follows: 
We start with a review of the theory of conformal superalgebras,
following mainly the lecture notes by Kac \cite{kac2} and 
also \cite{freecv}. In \sec{examples} we construct most of the
examples of conformal superalgebras used later on. Then in
\sec{vertex}-\sec{eva} we outline the theory of vertex algebras, using 
again \cite{kac2}. In \sec{lattice} we describe the 
construction of lattice vertex superalgebras. In \sec{tkk} we review
the Tits-Kantor-Koeher construction and use it to get the 
conformal algebra $\goth K$. Having done that we review the so-called
bozon--fermion correspondence, which is essentially the study of the lattice
vertex algebra $V_\Z$ corresponding to the lattice $\Z$. Sections
\sec{matrices}-\sec{fermion} are again taken from \cite{kac2}. 

In sections
\sec{v0}--\sec{repW} we explore the structure of $V_\Z$ as a module
over the conformal algebra $\goth W$ of differential operators on
a circle. This is the same as the module structure of the Lie algebra
$W_+ = \Bbbk\<\,t,p\,|\,\ad tp =1\,\>^{(-)}$ of differential operators
on a disk. We note that while the representation theory of the Lie
algebra $W= \Bbbk\<\,t,t,\inv,p\,|\,\ad tp =1\,\>^{(-)}$ of differential operators
on a circle, as well as that of the related vertex
algebra $\cal W_{1+\infty}$, has been extensively studied (see e.g. \cite{kr,fkrw}), 
the representation theory of $W_+$ seems to be mostly unknown. 

In \sec{rootsys} we give a rigorous formulation  of the problem and
introduce the necessary definitions. 
Then in \sec{rank1} we use the above results to study the subalgebras
of $V_\Lambda$ in case when $\op{rk} \Lambda = 1$. In \sec{rank2} we
proceed to the case when $\op{rk} \Lambda = 2$. This allows in
\sec{posdef}  to classify the root systems for the case
when the lattice $\Lambda$ is positive definite and in  \sec{finite}
to describe all finite root systems. Finally, in \sec{ears} we outline
the relation with the theory of EARS.  

\smallskip
Throughout this paper all spaces and algebras are over a ground field
$\Bbbk$ of characteristic 0.  

\subsection*{Acknowledgments}
I am very grateful to Igor Pak for teaching me elements of
combinatorics used in \sec{v0}--\sec{repW}. I also would like to thank
Stephen Berman, Yuly Billig, Jim Lepowsky, Weiqiang Wang and Efim
Zelmanov for helpful conversations concerning the present work.

%% file: conformal.tex
\subsection{Formal series and conformal algebras}
\label{sec:formser} 

Let $L=L\even \oplus L\odd$ be a Lie superalgebra. Consider the space of
formal power 
series $L[[z^{\pm1}]]$. We will write an element $\alpha\in
L[[z^{\pm1}]]$ in the form
$$
\alpha= \sum_{n\in\Z} \alpha(n)\,z^{-n-1}, \qquad \alpha(n)\in L. 
$$

Denote $L[[z^{\pm1}]]' = L\even[[z^{\pm1}]]\oplus
L\odd[[z^{\pm1}]]\subseteq L[[z^{\pm1}]]$.
The space $L[[z^{\pm1}]]$ is endowed with a derivation $D=d/dz$ and a
family of bilinear products $\ensquare n, \ n\in \Z_+$, given by 
\begin{equation}\label{fl:seriesprod}
\big(\alpha\ensquare n \beta\big)(m) = 
\sum_{i=0}^n \binom ni \ad{\alpha(n-i)}{\beta(m+i)}.
\end{equation}

We say that a pair of formal series $\alpha,\beta \in L[[z^{\pm1}]]$
are {\it local} if there is $N = N(\alpha, \beta)\in\Z_+$ such that 
$$
\sum_{i=0}^N (-1)^i \binom Ni \ad{\alpha(n-i)}{\beta(m+i)}=0
$$ 
for all $m,n\in\Z$. In particular we have $\alpha \ensquare n \beta =
0 $ for all $n\ge N$. 

The Dong's lemma \cite{kac2,li} states that if
$\alpha,\,\beta,\,\gamma \in L[[z^{\pm1}]]$ are three pairwise local 
formal series, then $\alpha \ensquare n \beta$
and $\gamma$ are local for all $n \in \Z_+$.

Let $\goth L\subset L[[z^{\pm1}]]'$ be a subspace of pairwise 
local formal series
closed under $D$ and under all products $\ensquare n$. Then $\goth L$ is a
{\it Lie conformal superalgebra}. 

Alternatively, we can define a Lie conformal superalgebra axiomatically
as a $\Bbbk[D]$-module $\goth L=\goth L\even \oplus\goth L\odd$ 
equipped with a family of products 
$\ensquare n, \ n\in \Z_+$ satisfying the following  axioms
(see e.g. \cite{kac2,freecv}): For any homogeneous $a,b,c \in \goth L$,
\begin{itemize}
\item[C1.](locality)
$a\ensquare n b = 0$ for $n\gg0$;
\item[C2.]
$(Da)\ensquare n b = -n a\ensquare{n-1} b$;
\item[C3.]
$D(a\ensquare n b) = (Da) \ensquare n b + a \ensquare n (Db)$;
\item[C4.](quasisymmetry) 
\begin{equation}\label{fl:qs}
a\ensquare n b = -(-1)^{p(a)p(b)} \sum_{i\ge0} (-1)^{n+i}\frac 1{i!} D^i (b\ensquare{n+i}a);
\end{equation}
\item[C5.](conformal Jacoby identity)
\begin{align}\label{fl:confjac}
(a \ensquare{n} b)&\ensquare m c =\\ 
&\sum_{i=0}^n (-1)^i\binom{n}{i} 
\biggl(a\ensquare{n-i} (b \ensquare{m+i} c) -(-1)^{p(a)p(b)}
b\ensquare{m+i} (a\ensquare{n-i} c)\biggr).\notag
\end{align}
\end{itemize}
Here $p(a)$ is the parity of $a$. 

One can prove that any subspace
in $L[[z^{\pm1}]]'$ of pairwise local series, closed under the products
\fl(seriesprod) and $D=d/dz$ satisfies all these axioms.

We will often use the notation $D^{(n)} = (-1)^n\frac 1{n!}\, D^n$.

\subsection{The coefficient algebra}
\label{sec:coeff}
Let $\goth U$ be a $\Bbbk[D]$-module. Its {\it space of coefficients} 
$U =\cff \goth U$ is constructed as follows. Consider the space 
$\goth U\otimes\Bbbk[t,t\inv]$, where $t$ is an independent
variable. We will write  
$a\otimes t^n = a(n)$ for $a\in \goth U$. Let 
$E=\op{Span}_\Bbbk \set{(Da)(n)+n\,a(n-1)}{a\in\goth U,\ n\in\Z}$.
Then let 
$$
U = \cff\goth U =\goth U\otimes\Bbbk[t,t\inv]/E.
$$
There is a homomorphism $\goth U\to U[[z^{\pm1}]]$  given by 
$a\mapsto \sum_n a(n)\,z^{-n-1}$. This homomorphism is the universal
one among all  the representations of $\goth U$ by formal series: if 
$\goth U\to U'[[z^{\pm1}]]$ is another $\Bbbk[D]$-ho\-mo\-mor\-phism, 
then there is a
unique homomorphism $U\to U'$ such that the diagram 
\begin{equation*}
\begin{array}{c}
U[[z^{\pm1}]] \xrightarrow{\hspace{15pt}} U'[[z^{\pm1}]] \\[-5pt]
\text{\Large $\nwarrow$} \hspace{15pt} \text{\Large $\nearrow$}\\[-10pt]
\goth U
\end{array}
\end{equation*}
commutes.

If $\goth U$ has a structure of a conformal algebra, then the space $U=\cff \goth U$ becomes a
``usual'' algebra. The product on $U$ is defined by 
$$
\ad{a(m)}{b(n)}= \sum_{i\ge0} \binom mi \big(a\ensquare i b\big)(m+n-i).
$$
The sum here makes sense due to the locality of $a$ and $b$. 
In this case $U$ is called the {\it coefficient algebra} of $\goth U$. It
still has the  universality property mentioned above.

Let $\goth L$ be a Lie conformal superalgebra and let $L=
\cff\goth L$ be its Lie superalgebra of coefficients. Then 
$L = L_-\oplus L_+$ is a direct sum of subalgebras 
$L_- = \op{Span}\set{a(n)}{a\in \goth L, \ n<0}$ and 
$L_+ = \op{Span}\set{a(n)}{a\in \goth L, \ n\ge0}$. 
The derivation $D$ is also  a derivation of $L$, acting by 
$D(a(n)) = -n a(n-1)$. We see that $L_-$ and $L_+$ are closed under
the action of $D$.

\subsection{Conformal modules}
\label{sec:cm}
Let as before $\goth L$ be a Lie conformal superalgebra and let $L=
\cff\goth L= L_-\oplus L_+$ be its Lie superalgebra of coefficients.
Denote by $\^L_+ = L_+ \oplus \Bbbk D$ the extension of $L_+$ by $D$.
A module  over $\goth L$ is by definition a $\^L_+$-module $\goth U$,
such that for any $u\in \goth U$ and $a\in \goth L$ we have $a(n)u=0$
for $n\gg0$. One can view $\goth U$ as a $\Bbbk[D]$-module such that
for any $a\in \goth L, \ u\in \goth U$ and $n\in \Z_+$ there is an
action $a\ensquare n u \in \goth U$, so that the semidirect product
$\goth L \ltimes \goth U$ becomes a Lie conformal superalgebra, and
$\goth U$ being its abelian ideal. 

The space of coefficients $U = \cff\goth U$ becomes a module over $L$
by 
$$
a(m)\, u(n) = \sum_{i\ge0} \binom mi \big(a\ensquare i u\big)(m+n-i).
$$

%% file: examples.tex
\subsection{Examples}\label{sec:examples}
\subsubsection{Affine algebras}
\label{sec:affine}
Let $\goth g$ be an arbitrary Lie superalgebra. Consider the corresponding
loop algebra  
$\~{\goth g} =\goth g\otimes\Bbbk[t,t\inv]$.
Now for any $a\in \goth g$, define 
$$
\~a = \sum_{n\in \Z} a t^n\, z^{-j-1} \in 
\~{\goth g}\,[[z, z\inv]].
$$

It is easy to see that any two $\~a, \~b$ are local with $N(\~a, \~b) = 1$
and 
$$
\~a\ensquare{0} \~b = \~{\ad ab}.
$$
By the Dong's lemma the series 
$\{\~a \ | \ a\in \mathfrak g\} \subset \~{\mathfrak g}\,[[z, z\inv]]$ 
generate a Lie conformal superalgebra $\goth G$. As a $\Bbbk[D]$-module,
$\goth G\cong \Bbbk[D]\otimes\goth g$.

In practice we are often interested in central extensions of loop
algebras. Assume that $\goth g$ has trivial odd part, and that 
$\goth g$ is equipped with an invariant bilinear form $(\,\cdot\,|\,\cdot\,)$.
Consider then  the Lie algebra 
$\op{Aff}(\goth g)=(\goth g\otimes\Bbbk[t,t\inv])\oplus \Bbbk \c$ with the brackets given by 
$$
\ad{a(m)}{b(n)} = \ad ab(m+n) + \delta_{m,-n}\, m\, (a|b)\,\c.
$$
The algebra $\op{Aff}(\goth g)$ is called the {\it affinization} of $\goth g$.
It is the coefficient algebra of a conformal algebra
$\goth{Aff}(\goth g)\subset G[[z,z\inv]]$ which is generated 
by the series $\~a= \sum_n a(n)\,z^{-n-1}$ for $a\in \goth g$ and
$\c=\c(-1)$ so that $D\c=0$ and
$$
\~a\ensquare{0} \~b = \~{\ad ab},\qquad
\~a\ensquare{1} \~b = (a|b)\,\c.
$$

In the case when $\goth g$ is an abelian Lie algebra, the
corresponding Affine algebra $\op{Aff}(\goth g)$ is a Heisenberg
algebra, and $\goth{Aff}(\goth g)$ is a Heisenberg 
conformal algebra. In the physics literature the series $\~a$ are
sometimes referred to as {\it bozons} in this case. 

\subsubsection{The Clifford  algebra}
\label{sec:clifford}
As  another example, take $\goth g$ to be a two-dimensional odd
linear space  spanned over $\Bbbk$ by $g_1$ and $g_{-1}$. Consider the
central extension $Cl=\goth g\otimes \Bbbk[t,t\inv]\oplus  \Bbbk \c$ of
the corresponding loop algebra with the brackets given by 
$$
\ad{g_\e(m)}{g_{-\e}(n)}= \delta_{m+n,-1}\c,
\quad\e=\pm1,
$$
the rest of the brackets are 0. We let $\c$ to be even.
The algebra $Cl$  is called the Clifford Lie superalgebra.
It is the coefficient algebra of the conformal Lie superalgebra 
$\goth{Cl}\subset Cl[[z,z\inv]]$ spanned over $\Bbbk[D]$ by $\gamma_\e =
\~g_\e,\ \e=\pm1$, and $\c=\c(-1)$ with the products given by 
$\gamma_{\e}\ensquare 0 \gamma_{-\e}=\c$  
(the rest of the products are 0). The series $\gamma_{\pm\e}$ are
sometimes called {\it fermions} by physicists.

The Clifford algebra $Cl$ is doubly graded: set $p(\gamma_\e(n))=\e$
 and $d(\gamma_\e(n))= -n-\frac12$ for $\e=\pm1, \ n\in\Z$, so we get 
$$
Cl=\bigoplus_{p\in\Z}Cl_p,\qquad
Cl_p=\bigoplus_{d\in\Z/2} Cl_{p,d}.
$$
The conformal Clifford algebra $\goth{Cl}$ is also doubly graded such that
$d(\gamma_\e)=\frac 12,\ p(\gamma_\e)=\e, \ p(\c)=d(\c)=0$ and 
$\goth{Cl}_{p,d}\ensquare n \goth{Cl}_{p',d'} \subset 
\goth{Cl}_{p+p', d+d'-n-1}$, \ 
$D\goth{Cl}_{p,d} \subset \goth{Cl}_{p,d+1}$.

\subsubsection{The Lie  algebra of differential operators}
\label{sec:diffops}
Another example of conformal algebras is obtained from the Lie algebra of
differential operators. Let $A = \Bbbk\<p,t^{\pm1}\,|\,\ad tp=1\,\>$ be
the (localization in $t$ of) the associative Weyl algebra. Let 
$W = A^{(-)}$ be the corresponding Lie algebra. It is the coefficient
algebra of a Lie conformal algebra $\goth W\subset W[[z,z\inv]]$ which
is spanned over $\Bbbk[D]$ by elements
$$
p_m = \sum_{n\in\Z}\frac 1{m!}\, p^m t^n\, z^{-n-1}\subset W[[z,z\inv]]
$$  
with the multiplication table
\begin{equation}\label{fl:confpp}
p_m\ensquare k p_n = 
\binom{m+n-k}m p_{m+n-k} -
(-1)^k\sum_{s=0}^{m-k}   \binom{m+n-k-s}n D^{(s)} p_{m+n-k-s}.
\end{equation}

The algebra $W$ has a unique central extension $\^W=W\oplus \Bbbk \c$,
defined by the 2-cocycle $\phi:W\times W\to \Bbbk$ given by 
\begin{equation}\label{fl:PhiW}
\phi(p^mt^k, p^nt^l) = \delta_{m+n,k+l}\,(-1)^m m!\,n!\binom k{m+n+1}.
\end{equation}
The algebra $\^W$ is usually referred to as $\cal W_{1+\infty}$, see
e.g. \cite{fkrw}. It is the coefficient algebra of a central extension
$\^{\goth W}=\goth W\oplus \Bbbk \c$ of the conformal Lie algebra
$\goth W$, defined by the conformal 2-cocycle \cite{bkv}
$\phi_k:\goth W\times \goth W \to \Bbbk\c, \ k\in\Z_+,$ given by 
\begin{equation}\label{fl:confphiW}
\phi_k\big(p_m,p_n\big) = \delta_{k,m+n+1}(-1)^m \c.
\end{equation}

The algebras $W$ and $\^W$ are graded by setting 
$\deg p = \deg t\inv =1,\ \deg t=-1,\ \deg \c =0$.   
The conformal algebras $\goth W$ and $\^{\goth W}$ inherit the
gradation from $W$ and $\^W$ respectively, so that we have 
$\deg(p_m) = m+1,\ \deg D=1,\ \deg \ensquare k = -k-1$.

\subsubsection{The  Virasoro conformal algebra}
\label{sec:virasoro}
Here are all non-zero products in $\^{\goth W}$ between $p_0$ and $p_1$:
\begin{gather*}
p_0\ensquare 1 p_0 = \c,\quad
p_0\ensquare 1 p_1 = p_0, \quad
p_0\ensquare 2 p_1 = \c,\\
p_1\ensquare 0 p_1 = Dp_1, \quad
p_1\ensquare 1 p_1 = 2p_1, \quad
p_1\ensquare 3 p_1 = -\c,\\
p_1\ensquare 0 p_0 = Dp_0, \quad
p_1\ensquare 1 p_0 = p_0, \quad
p_1\ensquare 2 p_0 = -\c.
\end{gather*}
The element $p_0\in\^{\goth W}$ generates a copy of the 
Heisenberg conformal algebra 
$\goth H=\goth{Aff}(\Bbbk)\subset \^{\goth W}$, introduced in \sec{affine}.
Its coefficient algebra $H=\cff\goth H$ is the affinization of the
one-dimensional trivial Lie algebra, so we have 
$$
\ad{p_0(m)}{p_0(n)} = \delta_{m,-n} m\, \c.
$$

The element $p_1\in\^{\goth W}$ generates 
the Virasoro conformal algebra $\goth{Vir}\subset   \^{\goth W}$. 
Its coefficient algebra $Vir = \cff\goth{Vir}$ is spanned by $p_1(m)$
for $m\in\Z$ and $\c$ with the brackets given by 
$$
\ad{p_1(m)}{p_1(n)} = (m-n)\,p_1(m+n-1) - \delta_{m+n,2} \binom m3\c.
$$
Together $p_0$ and $p_1$ span a semidirect product 
$\goth{Vir}\ltimes \goth H\subset  \^{\goth W}$.

Note that $\deg p_0 = 1$ and $\deg p_1 = 2$.

\subsubsection{$N=2$ simple Lie conformal superalgebra}
\label{sec:N=2}
A conformal algebra $\goth A$ is said to be of a {\it finite type} if it is a
finitely generated module over $\Bbbk[D]$. 
The algebras $\goth{Vir}, \goth{Cl}$ and $\goth{Aff}(\goth g)$ for
finite dimensional $\goth g$ defined above are of a finite type. 
All simple and semisimple Lie conformal superalgebras of finite type
are classified by Kac in \cite{kac98}, see also \cite{dk} for the
non-super case. 

Besides the algebras mentioned
above we will need in the sequel the following simple finite type Lie conformal
superalgebra, called the $N=2$ Lie conformal
superalgebra. It is spanned over $\Bbbk[D]$ by two odd elements 
$\gamma_{-1}, \gamma_{+1}$ and two even elements $v, \ h$. 
Elements $v$ and $h$ generate respectively the Virasoro and Heisenberg
Lie conformal algebra, so the even part of the $N=2$ superalgebra is
equal to $\goth{Vir}\ltimes \goth H$. The remaining non-zero
products between the generators are
\begin{gather*}
\gamma_{-1}\ensquare 0 \gamma_1 = v + \frac 12 Dh,\ \ 
\gamma_{-1}\ensquare 1 \gamma_1 = h,\ \
v\ensquare 0 \gamma_{\pm1} = D\gamma_{\pm1},\\ 
v\ensquare 1 \gamma_{\pm1} = \frac 32 \gamma_{\pm1},\ \ 
h\ensquare 0 \gamma_{\pm1} = \pm \gamma_{\pm1}.
\end{gather*}
Use \fl(qs) to get products in the other order.

%% file: vertex.tex
\subsection{Vertex algebras}
\label{sec:vertex}
In order to define vertex algebras, we need to consider the so-called
{\it fields} instead of formal power series. Let $V=V\even \oplus V\odd$ be
a vector superspace. Denote by $gl(V)$ the Lie
superalgebra of all $\Bbbk$-linear operators on $V$. Consider the space 
$\F(V) \subset gl(V) [[z, z^{-1}]]'$
of {\it fields} on $V$, given by
$$
\F(V) = \left\{ \left.
\sum_{n\in \mathbb Z}a(n)\,z^{-n-1}\ \right|
\forall v \in V, \  a(n) v = 0 \text{  for } n \gg 0
\right\}.
$$
For $\alpha(z) \in \F(V)$ denote 
$\alpha_-(z) = \sum_{n<0}\alpha(n)\,z^{-n-1}$,\ 
$\alpha_+(z) = \sum_{n\ge 0}\alpha(n)\,z^{-n-1}$.
Denote also  by  $ \1 = \1_{\F(V)} \in \F(V)$ the identity operator, 
such that $ \1 (-1) = \op{Id}_V$,  all other coefficients are 0.

In addition to the products $\ensquare n, \ n\in \Z_+$, defined by
\fl(seriesprod), the space of fields $\F(V)$ has products $\ensquare
n$  for  $n<0$. Define first $\ensquare{-1}$ by 
\begin{equation*}
\alpha(z) \ensquare{-1} \beta(z) =  \alpha_-(z) \beta(z) + \beta(z) \alpha_+(z).
\end{equation*}
Note that the products of fields here make sense, since for any 
$v \in V$ we have $\alpha(n)v = \beta(n)v = 0$ for $n \gg 0$.

Next, for any $n<0$ set 
\begin{equation*}
\alpha(z) \ensquare{n} \beta(z) = 
\frac1{(-n-1)!}\,  \big(D^{-n-1}\alpha(z)\big)\ensquare{-1} \beta(z) ,
\end{equation*}
where $D = \frac d{dz}$. 
It is easy to see that 
\begin{equation*}
\alpha \ensquare{-1} \1 = \alpha, \qquad \alpha \ensquare{-2} \1 = D \alpha, 
\qquad \1 \ensquare n \alpha = \delta_{-1,n} \alpha.
\end{equation*}

It also follows easily from definitions that $D$ is a
derivation of all these products:
\begin{equation*}
D(\alpha\ensquare n \beta) = D \alpha \ensquare n \beta + \alpha \ensquare n D \beta.
\end{equation*}

We have the following explicit formula for the products:
if $\big(\alpha \ensquare{n} \beta\big)(z)
 = \sum_m \big(\alpha \ensquare{n} \beta\big)(m)\, z^{-m-1}$, then 
\begin{equation*}
\begin{split}
\big(\alpha \ensquare{n} \beta\big)(m) 
& = \sum_{s\le n} (-1)^{s+n} \binom n{n-s} \alpha(s)\beta(m+n-s)\\
& -(-1)^{p(a)p(b)} \sum_{s\ge0} (-1)^{s+n} \binom n{s} \beta(m+n-s)\alpha(s).
\end{split}
\end{equation*}

The notion of locality introduced in \sec{formser} applies also to
fields without any changes. The Dong's lemma holds for fields
instead of formal power series as well.

A {\it vertex superalgebra} is a subspace of fields $\goth V \subset
\F(V)$ such that $\1 \in \goth V$, \ $\alpha \ensquare n \beta \in
\goth V$ for every $\alpha, \beta \in \goth V$ and $n \in \Z$, and 
every $\alpha, \beta \in \goth V$ are local to each other. For an
axiomatic definition of vertex superalgebras, equivalent to the above
description, we refer the reader to \cite{bor,dl,fhl,flm,kac2}. 

For a vertex superalgebra $\goth V$ one can consider the {\it left regular
action map} $Y:\goth V \to \F(\goth V)$ given by 
$Y(a)(z) = \sum_{n\in\Z} (a\ensquare n \cdot\,)\, z^{-n-1}$. One of the
main properties of vertex superalgebras is that $Y$ is a vertex superalgebra
homomorphism, in particular 
$Y(a\encirc n b) = Y(a) \encirc n Y(b)$, which is equivalent to the
following generalization of the conformal Jacoby identity \fl(confjac):
\begin{equation}\label{fl:verjac}
\begin{split}
\bigl(a \ensquare n b\bigr) \ensquare m c &= 
 \sum_{s\le n} (-1)^{s+n} \binom n{n-s} a\ensquare s \bigl(b
\ensquare{m+n-s} c\bigr)\\
& -(-1)^{p(a)p(b)} \sum_{s\ge0} (-1)^{s+n} \binom n{s}
b\ensquare{m+n-s}\bigl(a\ensquare s c\bigr),
\end{split}
\end{equation}
for all $m,n \in \Z$.

We note that vertex superalgebras are in particular conformal
superalgebras. Moreover, it can be shown that the quasisymmetry
identity \fl(qs) holds in vertex superalgebras for all integer $n$.

\subsection{Enveloping vertex algebras of a conformal algebra}
\label{sec:eva}
Let again $\goth L$ be a conformal superalgebra and 
$L=\op{Coeff}\goth L=L_-\oplus L_+$ be its 
coefficient Lie superalgebra. Denote by $\widehat L = L\oplus \Bbbk D$ the
extension of $L$ by the derivation $D$, see \sec{coeff}.
Consider a homomorphism of conformal superalgebras
$\rho:\goth L\to \goth V$ of $\goth L$ into a
vertex superalgebra $\goth V$.
If $\goth V$ is generated as a vertex
superalgebra by $\rho(\goth L)$ then we call it an {\it enveloping vertex
superalgebra} of $\goth L$. 

An $L$-module (or $\widehat L$-module) $U$ is called a {\it highest
weight} module if  it is generated as a module over $L$  by a single
element $u \in U$  such that  
$L_+ u = Du=0$. In this case $u$ is called a {\it highest weight vector}.

It is well-known \cite{kac2,primc,freecv} that 
the enveloping vertex superalgebra $\goth V$ has the structure of a
highest weight module over $\widehat L=\cff\goth L\oplus \Bbbk D$ 
with the highest weight vector $\1$ defined by
$a(n)v = \rho(a)\ensquare n v$.
Ideals of $\goth V$ are $\widehat L$-submodules and if  $\rho_1:\goth L\to
\goth V_1$ and $\rho_2:\goth L\to \goth V_2$ are two enveloping
vertex algebras of $\goth L$ and $\psi:\goth V_1 \to \goth V_2$ is a
vertex algebra homomorphism  such that $\rho_1 \psi = \rho_2$ then 
$\psi$ is an $\widehat L$-module homomorphism. 
Conversely, any highest weight module $V$ over $\widehat L$ with the highest
weight vector $\1$ has a structure of enveloping vertex algebra of
$\goth L$ with the map $\rho:\goth L\to V$ given by 
$\rho(a) = a(-1)\1$. In this case we have 
$a(n)v = \rho(a)\ensquare n v$ for $a\in\goth L$,
submodules of $V$ are vertex ideals and if $\rho:V_1\to V_2$ is a
homomorphism of two 
highest weight $L$-modules such that $\rho(\1)=\1$ then $\rho$ is a
vertex algebra homomorphism.

We also mention the notion of universal (or Verma)
highest weight module over $L$.  It is defined by
$V(L) = \op{Ind}_{L_+}^L \Bbbk\1 = U(L)\otimes_{U(L_+)}\Bbbk\1$. 
The action of the derivation $D$ on $L$ can be naturally extended to 
the action on $V(L)$, so $V(L)$ becomes an $\widehat L$-module.
Verma module is universal in the sense that for any other highest
weight module $U$ with highest weight vector $u$ there is unique
homomorphism $V(L)\to U$ such that $\1\mapsto u$. So the theorem
implies that the enveloping vertex algebra corresponding to the Verma
module $V(L)$ is universal in the obvious sense. It is called the 
{\it universal enveloping vertex algebra} of $\goth L$.

%% file: lattice.tex
\subsection{Lattice vertex algebras}
\label{sec:lattice}
In this section we construct a very important example of vertex
superalgebras, called lattice vertex superalgebras.  We mostly follow 
\cite{kac2}, see also \cite{dong,fk,flm}. The
Frenkel-Lepowsky-Meurman Moonshine vertex algebra $V^\natural$, for
example, is closely related to the vertex algebra corresponding to the
Leech lattice --- a simple unimodular lattice of rank 24. Also, the
lattice vertex algebras play a very important
role in physics. 

We start from  an integer lattice $\Lambda$ of rank
$\ell$. It comes with
an integer-valued bilinear form $(\,\cdot\,|\,\cdot\,)$. Let 
$\goth h = \Lambda \otimes_\Z \Bbbk$ and we extend the form to $\goth h$. 
Let $\nu:\goth h \to \goth h^*$ be the usual identification of 
$\goth h$ with $\goth h^*$ by means of this form. 

Let $\goth H = \goth{Aff}(\goth h)$
be the Heisenberg conformal algebra corresponding to the space
$\goth h$ (see \sec{affine}), and let $H=\op{Aff}(\goth h) = \cff \goth H$
be its coefficient Heisenberg algebra.
Take some  $\beta \in \Lambda$ and let $V_\beta$ be the canonical relation
representation of $H$, that is, a highest weight irreducible 
$H$-module generated by the highest weight vector $v_\beta$ such that 
$h(0)=(h|\beta)\,\op{Id}$,\ $\c=\op{Id}$.

It follows from
\sec{eva} that $V_0$ is
an enveloping vertex algebra of $\goth H,\ v_0 = \1$. It can be shown that  
$V_\beta$ is a module over  the vertex algebra $V_0$. 
We define $V_\Lambda=\bigoplus_{\beta\in\Lambda} V_\beta = V_0\otimes\Bbbk[\Lambda]$. 

Let $\e:\Lambda\times \Lambda \to \{\pm 1\}$ be a
bimultiplicative map such that 
\begin{equation}\label{fl:eps}
\e(\alpha,\beta) = (-1)^{(\alpha|\alpha)(\beta|\beta)} 
(-1)^{(\alpha|\beta)}\e(\beta,\alpha).
\end{equation}
for any $\alpha,\beta\in\Lambda$. We remark that it is enough to check 
the identity \fl(eps) only when $\alpha$ and $\beta$ belong to some 
$\Z$-basis of $\Lambda$; then \fl(eps) will follow for general
$\alpha, \beta$ by bimultiplicativity.

Let 
$$
1\longrightarrow \{\pm1\} \longrightarrow \^\Lambda \longrightarrow \Lambda
\longrightarrow 1
$$
be the extension of $\Lambda$ corresponding to the cocycle $\e$. Let 
$e:\Lambda\to \^\Lambda$ be a section of this extension. The extended lattice
$\^\Lambda$ acts on the group algebra $\Bbbk[\Lambda]$ of $\Lambda$ by 
$e(\alpha)e^\beta = \e(\alpha,\beta) e^{\alpha+\beta}$.

For $m\in \Z_+$ let 
$\cal P(m) = \set{k=(k_1,k_2,\ldots)}{k_i\ge 0, \ \sum_{i\ge1} i k_i=m}$ 
be the set of partitions of $m$.

The main result \cite{dong,dl,flm,kac2} is that there is a unique
vertex superalgebra 
structure on $V=V_\Lambda$ such that $a \ensquare n v = a(n)v$ for any 
$a\in \goth H \subset V_0$ and $v\in V$. The products are defined by 
\begin{equation}\label{fl:vanvb}
v_\alpha \ensquare n v_\beta = 
\e(\alpha,\beta)\sum_{k\in \cal P(-(\alpha|\beta)-n-1)} \ 
\prod_{j\ge1} 
\(\frac{\alpha(-j)}{j!}\)^{k_j} v_{\alpha+\beta}.
\end{equation}
In particular, $v_\alpha \ensquare n v_\beta = 0$ if 
$n \ge -(\alpha|\beta)$ and 
$v_\alpha \ensquare{-(\alpha|\beta)-1} v_\beta = 
\e(\alpha,\beta)\, v_{\alpha+\beta}$, \ 
$v_\alpha \ensquare{-(\alpha|\beta)-2} v_\beta = 
\e(\alpha,\beta)\, \alpha(-1)v_{\alpha+\beta}$. Note that $V_\alpha
\ensquare n V_\beta \subset V_{\alpha+\beta}$.

The even and odd parts of $V$ are
$$
V\even = 
\bigoplus_{\substack{\alpha\in \Lambda:\\ (\alpha|\alpha)\in 2\Z}} V_\alpha, 
\qquad
V\odd = 
\bigoplus_{\substack{\alpha\in \Lambda:\\ (\alpha|\alpha)\in 2\Z+1}}
V_\alpha.
$$
The vertex algebra $V$ is simple if the form
$(\,\cdot\,|\,\cdot\,)$ is non-degenerate.

Under the left regular action map $Y:V \to \F(V)$ 
the elements $v_\alpha$ are mapped to the so-called vertex operators
$$
Y(v_\alpha)=\Gamma_\alpha(z) = e(\alpha) z^{\alpha(0)} 
E_-(\alpha, z) E_+(\alpha, z),
$$
where 
$$
E_\pm(\alpha, z) = \exp \sum_{n\gtrless 0}
- \frac{\alpha(n)}n z^{-n}\in\F(V),
$$
and the field $z^{\alpha(0)}\in \F(V_\Z)$ is defined by 
$z^{\alpha(0)}\big|\raisebox{-4pt}{$V_\mu$} = 
\sum_{n\in \Z}\delta_{n, (\alpha|\mu)} \, z^n$.
We also have
$$
\ad{h(n)}{e(\alpha)} = \delta_{n,0}\, (\alpha|h)\, e(\alpha).
$$

Besides the grading by the lattice $\Lambda$, the vertex superalgebra 
$V$ has another grading by the group $\frac 12\Z$, so that 
$\deg v_\alpha = \frac12 (\alpha|\alpha),\ \deg a = 1$ for every
$a\in \goth H,\ \deg \ensquare n = -n-1$ and $\deg D = 1$. We have
decomposition 
$$
V_\beta = 
\bigoplus_{d\in \frac{(\beta|\beta)}2 + \Z_+} V_{\beta,d}.
$$

Let $(\alpha_1,\ldots,\alpha_\ell)$ and $(\beta_1,\ldots,\beta_\ell)$
be dual bases of $\goth h$, i.e. such that
$(\alpha_i|\beta_j)=\delta_{ij}$. Then the element 
$\omega = \frac 12\sum_{i=1}^\ell \alpha_i \ensquare{-1} \beta_i \in V_0$
generates a copy of the Virasoro Lie conformal algebra $\goth{Vir}$,
defined in \sec{virasoro}, such that $\omega\ensquare 0 v = DV$ for 
all $v\in V$, \  $\omega \ensquare 1 v = (\deg v) v$ for all
homogeneous $v\in V$, \ $\omega \ensquare 2\omega =0$ 
and $\omega \ensquare 3 \omega = \frac 12 \1$.

We will identify the Heisenberg conformal algebra $\goth H$ with
its image in $V_0$ under the map $\~h\mapsto h(-1)\1$ for 
$h\in\goth h$,\ $\c\mapsto \1$. 

We remark that the vertex superalgebra structure of $V_\Lambda$ is very
explicit. A basis of $V_\Lambda$ is given by all expressions of the
form
$$
\alpha_1(n_1)\alpha_2(n_2)\ldots\alpha_l(n_l)v_\beta,\qquad 
\alpha_i, \beta \in \Lambda,\quad
0>n_i\in\Z, 
$$ 
and the products of these elements are easily calculated using the
formula \fl(vanvb), the following identities in $\cff V_\Lambda$:
$$
\ad{\alpha(m)}{v_\beta(n)}=(\alpha|\beta)\,v_\beta(m+n),\quad
\ad{\alpha(m)}{\beta(n)}=(\alpha|\beta)\,m\,\delta_{m,-n}
$$
for $\alpha,\beta \in \Lambda, m,n\in\Z$,
 and the identities \fl(qs) and \fl(verjac) of vertex superalgebras.

%% file: tkk.tex
\subsection{The Tits-Kantor-Koeher construction}
Let $L=L_{-1}\oplus L_0 \oplus L_1$ be a three-graded Lie algebra,
such that $\ad{L_i}{L_j} \subset L_{i+j}$ whenever
$i,j,i+j\in\{-1,0,1\}$ and $\ad{L_1}{L_1}=\ad{L_{-1}}{L_{-1}}=0$. 
Assume that $L_0 = \ad{L_{-1}}{L_1}$ and 
$L_0\cap Z(L)=0$, where $Z(L)$ is the center of $L$.
Consider a pair of trilinear maps
$$
\f_+:L_1\times L_{-1}\times L_1 \to L_1, \qquad
\f_-:L_{-1}\times L_1\times L_{-1} \to L_{-1}, 
$$
given by $\f_\pm(a, b, c) = \frac 12 \,\ad{\ad ab}c$. The tuple 
$J=\{\,(L_{-1}, L_1),\ \f_\pm\,\}$ is a so-called a {\it Jordan
pair}. All Jordan pairs can be obtained in this way. 
In fact one can define Jordan pairs formally by imposing certain axioms on the maps 
$\f_{\pm1}$. From an abstractly
defined Jordan pair $J=\{\,(L_{-1}, L_1),\ \f_\pm\,\}$ one  
can construct a three-graded Lie algebra 
$L(J) = L_{-1}\oplus L_0\oplus L_1$, where $L_0=\cal D(J)$ 
is the Lie algebra of inner
derivations of $J$. This is is known as the
Tits-Kantor-Koeher construction, see \cite{tits,kantor,koeher}. 
A Jordan pair $J$ is simple if and only if the TKK Lie algebra $L(J)$
is simple.

A derivation $d = (d_-, d_+)$ of  a
Jordan pair $J=\{\,(L_{-1}, L_1),\ \f_{\pm1}\,\}$ is a pair of linear
maps $d_-:L_-\to L_-, \ d_+:L_+\to L_+$, such that  
$$
d_\pm\big(\f_\pm(a,b,c)\big) = 
\phi_\pm\big(d_\pm(a),b,c\big) + \f_\pm\big(a,d_\mp(b),c\big) +
\f_\pm\big(a,b,d_\pm(c)\big ).
$$
As it is the case for other algebraic structures, the set of all
derivations of $J$ is a Lie algebra under the usual commutator
operation. 
Let $a\in L_{-1}, \ b\in L_1$. It turns out that 
$d_{ab}=\big(\,\f_-(a,b,\ \cdot\,),\, \f_+(b,a,\ \cdot\,)\big)$
is a derivation of $J$, called an inner derivation. The Lie algebra
$L_0$ can be identified with the set 
$\cal D(J) = \set{d_{ab}}{a\in L_{-1},\ b\in L_1}$ 
of all inner derivations of $J$ by $d_{ab} = \frac 12 \ad ab$.


There is a very important case when $L_{-1}=L_1$ and  $\f_+ =
\f_-$. Then $J$ is called {\it a Jordan triple system}. In terms of
the three-graded Lie algebra $L=L(J)$ this means that there is an involution 
$\sigma:L\to L$ such that $\sigma(L_\e)=L_{-\e},\ \e=\pm1,$ 
is the identification of $L_{-1}$ and $L_1$. All Jordan pairs we deal
with in this paper are in fact Jordan triple systems.

\subsection{Example: associative algebras}
Here we consider some important examples of Jordan triple systems. Let 
$A$ be an associative algebra.  We define a triple operation
$\f:A\times A\times A\to A$ by $\f(a,b,c) = \frac 12
(abc+cba)$. The TKK Lie algebra $L(A)$ can be identified with a
subalgebra of the Lie algebra $gl_2(A)$ of $2\times2$ matrices over $A$
modulo the center:
$$
L(A) =\op{Span}_\Bbbk \left\{\left. 
\begin{pmatrix}
-dc& b\\
a&cd
\end{pmatrix}\in gl_2(A)/Z\big(gl_2(A)\big)\ 
\right|\ 
a,b,c,d \in A
\right\}.
$$
Here $L(A)_{-1}$ consists of lower triagular matrices, $L(A)_1$
consists of upper triangular matrices and $L(A)_0$ is the space of
all diagonal matrices in $L(A)$.
Quite often it happens that $L(A) = gl_2(A)/Z\big(gl_2(A)\big)$. 

Now assume that there is an involution or an anti-involution
$\tau:A\to A$ on $A$. Then both the set $A^\tau$ of $\tau$-stable 
elements and the set $A^{-\tau}$ of those elements that change sign
under the action of $\tau$ are closed under the triple operation. The
corresponding TKK Lie algebra $L(A^{\pm\tau})$ can be still
represented by $2\times2$ matrices over $A$.
Consider the case when $\tau:A\to A$ is an anti-involution. Then 
$$
L(A^{\pm\tau}) =\left\{\left. 
\begin{pmatrix}
-\tau(x)& b\\
a&x
\end{pmatrix}\ 
\right|\ 
a,b \in A^{\pm\tau}, \ x\in \cal D\subset A^{(-)}
\right\},
$$
where $\cal D \subset A^{(-)}$ 
is a Lie subalgebra of $A^{(-)}$
generated by all elements of the form $ab$ for $a,b \in
A^{\pm\tau}$. 
We see that  $\cal D=L(A^{\pm\tau})_0$ is precisely the Lie algebra of inner
derivations of the Jordan pair $(A^{\pm\tau}, A^{\pm\tau})$. 
It acts on $A^{\pm\tau}$ by $x.a = xa + a \tau(x)$ where
$x\in \cal D,\ a\in  A^{\pm\tau}$. 
The involution $\sigma : L(A^{\pm\tau}) \to L(A^{\pm\tau})$ acts by 
$\bigl(\begin{smallmatrix}-\tau(x)& b\\a&x\end{smallmatrix}\bigr)
\mapsto 
\bigl(\begin{smallmatrix}x&a \\b&-\tau(x)\end{smallmatrix}\bigr)$,
therefore,
$\sigma\big|\raisebox{-4pt}{$\cal D$} = -\tau$.

\subsection{The conformal algebra $\^{\goth K}$}
\label{sec:weyl}
Now let
$A = \Bbbk\<p,t,t\inv\,|\,\ad tp=1\,\>$ be
the localized Weyl algebra, the one we have dealt with in \sec{diffops}. It
has an anti-involution $\tau:A\to A$ defined by $\tau(t)=t, \
\tau(p)=-p$.  The space $J=A^{-\tau}$ is a Jordan triple subsystem
of $A$. It is easy to see that $J$ is simple. Let $K=L(J) = 
J \oplus \cal D(J) \oplus J$ be the
TKK Lie algebra corresponding to $J$. 
\begin{lemma}\label{lem:L0}\sl
The Lie algebra $K$ is
the coefficient algebra of a simple conformal algebra 
$\goth K$. 
\end{lemma}
\begin{proof}
First we construct the $\Bbbk[D]$-module $\goth J$, such that 
$J = \cff \goth J$. 
The anti-involution $\tau$ acts also on the $\Bbbk[D]$-module
$\goth A\subset A[[z,z\inv]]$, generated by the series 
$p_m = \sum_n p_m(n)\,z^{-n-1} \in  A[[z,z\inv]]$, see \sec{diffops}.
(In fact $\goth A$ has a structure of an {\it associative} conformal
algebra, see \cite{kac_fd,freecv}).
We let $\goth J=\goth A^{-\tau}$ to be the
$\Bbbk[D]$-submodule of $\goth A$ consisting of those elements
of $\goth A$ which change sign under the action of $\tau$. We have
$J=\cff\goth J$.

The Weyl algebra $A$ is spanned by the coefficients 
$p_m(n) = \frac 1{m!}\, p^m t^n,\ m\in\Z_+, \ n\in\Z$, see
\sec{diffops}. Therefore the space 
$J$ is spanned by the elements
$$
u_m(n) = \tau(p_m(n)) - p_m(n) = 
-\frac 1{m!}\,\, p^m t^n + 
(-1)^m \sum_{i\ge0} \binom ni \frac 1{(m-i)!}\,\, p^{m-i}t^{n-i}
$$
for all $m\in\Z_+, \ n\in\Z$. Then the series
$$
u_m = \sum_{n\in\Z} u_m(n)\,z^{-n-1} = 
-p_m+(-1)^m \sum_{i=0}^m D^{(i)} p_{m-i}\in \goth J
$$
span $\goth J$ over $\Bbbk[D]$. 
Since $u_m(n)+\tau(u_m(n))=0$, we have 
$$
u_m+(-1)^m \sum_{i\ge0} D^{(i)}u_{m-i}=0,
$$
therefore, if $m$ is even, then we get
\begin{equation}\label{fl:udepend}
u_m = - \frac 12 \sum_{i\ge 1}  D^{(i)}u_{m-i}.
\end{equation}
The remaining series $\set{u_m}{m\in 2\Z_++1}$ form a basis of 
$\goth J$ over $\Bbbk[D]$. 

Some simple calculations show that the  algebra $\cal D(J)$ of
inner derivations of  $J=A^{-\tau}$ is equal to
the whole $W=A^{(-)}$. So $J$ is a module over
$W$, where the action is given by 
$x.a = xa + a \tau(x)$ for $x\in W$ and $a\in J$. This action
induces the following action of $\goth W$ on $\goth J$:
\begin{equation}\label{fl:confpu}
p_m(k)u_n = \binom{m+n-k}m u_{m+n-k}
-\delta_{k,0}(-1)^n\sum_{i\ge0}  \binom {m+n-i}{m}\,D^{(i)}u_{m+n-i}
\end{equation}
for all $m,n\in\Z_+$.

So we obtain  TKK conformal algebra 
$\goth K = \goth K_{-1}\oplus \goth K_0 \oplus \goth K_1
\subset K[[z,z\inv]]$, where a $\Bbbk[D]$-basis of 
$\goth K_0 = \goth W$ is given by $p_m \in W[[z,z\inv]]$ for
$m\in\Z_+$, a $\Bbbk[D]$-basis of 
$\goth K_{-1} = \goth J$ is given by
$u_m \in J[[z,z\inv]]=K_{-1}[[z,z\inv]]$ for $m\in 2\Z_++1$, and 
$\goth K_1 = \goth J$ is spaned  over $\Bbbk[D]$ by the basis
$\sigma(u_m) \in K_1[[z,z\inv]]$ for $m\in 2\Z_++1$.
Here $\sigma:\goth K\to \goth K$ is the involution identifying 
$\goth K_1$ with $\goth K_{-1}$. Since the coefficients of these
series span $K$ and are linearly independent, we have 
$K = \cff\goth K$.

The formula \fl(confpu) shows that 
$p_m$ and $u_n$ in $K[[z,z\inv]]$ are local, hence so are 
$p_m$ and $\sigma(u_n)$. The series $u_m$ and $\sigma(u_n)$ are local
as well because the product $K_{-1}\times K_1 \to K_0$ is just 
the the associative product in  $A$ if we identify 
$K_{-1} = K_1 = A^{-\tau}\subset A$ and $K_0 = A$ as linear spaces. 
\end{proof}

Here is the multiplication table in $\goth K$. The products of
$p_m$ and $p_n$ are given by \fl(confpp), the products $p_m\ensquare k u_n =
p_m(k)u_n$ are given by \fl(confpu), 
\begin{equation}
\label{fl:confpsu}
\begin{split}
p_m\ensquare k \sigma(u_n) & = (-1)^{m+1} \binom {m+n}m \sigma(u_{m+n-k})\\
&+ (-1)^{m+n}\sum_{i=0}^n \binom{m+n-k-i}{m-k}\,
D^{(i)}\sigma(u_{m+n-k-i}),
\end{split}
\end{equation}
\begin{equation}\label{fl:confsuu}
\begin{split}
u_m\ensquare k \sigma(u_n) & =
\Biggl(\binom{m+n-k}m - (-1)^m \binom {m+n}m\Biggr) p_{m+n-k} \\
&-(-1)^n \sum_{j,l}(-1)^l \binom kl 
\Biggl(\binom{j-k}m - (-1)^m \binom {j-l}m\Biggr) D^{(m+n-j)}p_{j-k}.
\end{split}
\end{equation}

\begin{remark}
We see that $\goth J$ has all the rights to be
called a conformal Jordan triple system. We could have defined the
conformal triple operation on $\goth J$, this would be a
family of trilinear maps, depending on two integer parameters. 
\end{remark}

The algebras $K$ and $\goth K$ have central extensions 
$\^{K}=K \oplus \Bbbk\c$ and 
$\^{\goth K}= \goth K \oplus \Bbbk\c$ respectively, defined by  
2-cocycles $\phi:K\times K\to \Bbbk\c$ and 
$\phi_k:\goth K\times \goth K\to \Bbbk\c$, \ $k\in\Z_+$, given by 
\begin{equation}\label{fl:Phieven}
\begin{split}
\phi\big(u_m(k), \sigma(u_n(l))\big) &= 
\delta_{k+l,m+n}\,(-1)^m \binom k{m+n+1}\Biggl(\binom{m+n}m -1 \Biggr),\\
\phi_k\big( u_m, \sigma(u_n)\big) &=
\delta_{k,m+n+1}\, (-1)^m \Biggl(\binom{m+n}m -1 \Biggr).
\end{split}
\end{equation}

\begin{remark}
It is possible to show that conformal algebra $\goth K$
(and even its subalgebra $\goth W \ltimes \goth J$) is not embeddable
into an associative conformal algebra, though it is finitely generated
and has linear locality function, see \cite{universal}. It is proved
in \cite{universal} that the linearity of locality function is a
{\it necessary} condition for embedding of a finitely generated 
conformal algebra into an associative conformal algebra. 
\end{remark}

%% file: fermion.tex
\subsection{Lie algebras of matrices}
\label{sec:matrices}
Let $gl_\infty$ be the Lie algebra of infinite matrices which have
only finitely many non-zero entries. Denote by $E_{ij}$ the elementary
matrix that has the only non-zero entry at the $i$th row and $j$th
column. Then we have
\begin{equation*}
\ad{E_{i,j}}{E_{k,l}} = \delta_{jk}E_{il}-\delta_{il}E_{kj}.
\end{equation*}

Let $M$ be the Lie algebra of infinite matrices that have only finitely many
non-zero diagonals. 
Both algebras $gl_\infty$ and $M$ are graded by setting the degree of
$E_{ij}$ equal to $j-i$.

It is well-known \cite{kp81} that $gl_\infty$ and $M$ have unique central
extensions $\^{gl}_\infty = gl_\infty\oplus \Bbbk \c \subset \^M=M\oplus \Bbbk \c$ 
defined by the
2-cocycle $\phi(A,B) = \op{tr}\big(\ad AJ B\big)$, where 
$J=\sum_{i<0} E_{ii}\in M$. We set $\deg \c=0$.
The values of $\phi$ on the elementary
matrices are given by
\begin{equation}\label{fl:PhiM}
\phi(E_{i,j},E_{k,l}) = \begin{cases}
\delta_{il}\delta_{jk}&\text{if}\quad j<0\ \text{and}\ i\ge0,\\
-\delta_{il}\delta_{jk}&\text{if}\quad i<0\ \text{and}\ j\ge0,\\
0&\text{otherwise}.
\end{cases}
\end{equation}

The associative Weyl
algebra $A= \Bbbk\<t,t\inv,p\,|\,\ad tp=1\>$ is embedded into  the
associative algebra of infinite matrices by 
$$
t\mapsto \sum_{i\in\Z} E_{i,i-1}, \quad 
t\inv\mapsto \sum_{i\in\Z} E_{i,i+1},\quad 
p\mapsto -\sum_{i\in\Z}(i+1)\, E_{i,i+1},
$$ 
hence the Lie algebra $W=A^{(-)}$, defined in \sec{diffops}, is
embedded into the Lie algebra $M$. 

\begin{lemma}\label{lem:restr}\sl
The restriction of the 2-cocycle
$\phi$ on $W$ precisely coincides with the 2-cocycle on $W$ given by
\fl(PhiW). 
\end{lemma}
\begin{proof}
Express the linear generators of $W$ in terms of elementary matrices
\begin{equation}\label{fl:pmn}
p_m(n) = \frac 1{m!}\,p^m t^n = (-1)^m\sum_{i\in\Z}\binom{i+m}m E_{i,i+m-n}.
\end{equation}
Then using the formula \fl(PhiM) we obtain by some calculations that 
$$
\phi\big(p_m(k),p_n(l)\big) = \delta_{m+n,k+l}\,(-1)^m \binom k{m+n+1}.
$$
\end{proof}

Therefore the central extension $\^W$ of $W$ is embedded
into the Lie algebra $\^M$.

\subsection{The Clifford vertex superalgebra}
\label{sec:fermion}
Let $\goth{Cl}$ be the Clifford conformal superalgebra, defined in
\sec{clifford}, and let $Cl=\cff\goth{Cl}$ be its coefficient Lie
algebra. Let $U = U(Cl)/(\c=1)$ be the quotient of the universal enveloping
algebra of $Cl$ over the ideal generated by the relation $\c=1$. 

The following lemma is proved by a straightforward computation, see
e.g.  \cite{kac1}.

\begin{lemma}\label{lem:M0}\sl
Let $e_{ij} = \gamma_{-1}(i)\gamma_1(-j-1)\in U$ for $i,j\in\Z$ and let 
$$
\hat e_{ij} =  
\begin{cases}
-\gamma_1(-j-1) \gamma_{-1}(i)&\text{if}\quad 
i=j\ge0,\\
e_{ij}&\text{otherwise.}\quad 
\end{cases}
$$
Then the mapping $E_{ij}\mapsto e_{ij}$ defines an
embedding of the Lie algebra
$gl_\infty$ of infinite matrices  (see \sec{matrices}) into $U$, 
and the map $E_{ij}\mapsto \hat e_{ij},\ c\mapsto 1$ defines an
embedding of the Lie algebra $\^{gl}_\infty$ into $U$.
\end{lemma}

Note that $d(\hat e_{ij})=j-i=\deg E_{ij}, \ p(\hat e_{ij})=0$, see
\sec{clifford} for notations.

Now consider the universal highest weight module $V$ over $Cl$, see 
\sec{eva}. By
definition $V$ is generated over $Cl$ by a single element $\1$ such
that $\c\1=\1, \ Cl_+\1=0$ and $V = U(Cl)\otimes_{U(Cl_+)\oplus \Bbbk \c}\Bbbk\1$.
As a linear space $V$ can be identified with the Grassman algebra 
$\Bbbk[\gamma_\e(n)\ |\ \e=\pm1,\ n<0]$, on which $\c$ acts as
identity, $\gamma_\e(n)$  for $n<0$ acts by multiplication on the
corresponding variable, and if $n\ge 0$ then
$\gamma_\e(n)$ acts as an odd derivation. It follows that $V$ is an
irreducible $Cl$-module.

The module $V$ inherits the
double grading from $Cl$, so we have
$$
V=\bigoplus_{p\in\Z} V_p,\qquad
V_p=\bigoplus_{d\in\Z/2} V_{p,d}.
$$
It is easy to see that if $V_{p,d}\neq0$ then $d-\frac p2\in\Z$ and 
$d\ge \frac{p^2}2$. Indeed, let 
\begin{gather*}
w=\gamma_{-1}(n_1)\wedge\gamma_{-1}(n_2)\wedge\cdots\wedge\gamma_{-1}(n_k)
\wedge\gamma_1(m_1)\wedge\gamma_1(m_2)\wedge\cdots\wedge\gamma_1(m_l)\in V,\\
n_1<n_2<\ldots<n_k<0, \qquad m_1<m_2<\ldots<m_l<0,
\end{gather*}
be such that $p(w)=p>0$. Then $l\ge p$ and $d(w)\ge
-m_1-m_2-\ldots-m_l-\frac l2 \ge \frac{p^2}2$.

\begin{remark}
There is a alternative construction of $V$ using semi-infinite wedge
products, see e.g. \cite[chapter 14]{kac1}. It implies that in fact
$\dim_\Bbbk V_{p,d} = P\big(d-\frac{p^2}2\big)$, where $P(n)$ is the classical
partition function.
\end{remark}

The module $V$ has the structure of enveloping vertex superalgebra of
$\goth{Cl}$ (see \sec{eva}) such that the embedding $\goth{Cl}\to V$ is
given by $a\mapsto a(-1)\1$. We will identify $\goth{Cl}$ with its image
in $V$. We have 
$$
V\even = \bigoplus_{p\in2\Z} V_p,
\qquad
V\odd = \bigoplus_{p\in2\Z+1} V_p.
$$

A certain completion $\overline U$ of the algebra $U$ acts on the
vertex algebra
$V$ --- some infinite sums of the elements of $U$ make sense as
operators on $V$. In particular the closure of the algebra 
$\^{gl}_\infty\subset U$ spanned by $\hat{e}_{ij}$ (see \lem{M0}) is the algebra 
$\^M\subset \overline U$. It follows  also that the algebra $\^W\subset \^M$
acts on $V$ and there is a map $\^{\goth W}\to V$ given by 
$a\mapsto a(-1)\1$. The following lemma describes the image of
$\^{\goth W}\subset V$, see e.g. \cite{kac2}.

\begin{lemma}\label{lem:WinV}\sl
The mapping $p_m\mapsto \gamma_{-1}\ensquare{-m-1}\gamma_1$ defines an
embedding of the conformal algebra $\^{\goth W}$ into $V_0\subset V$.
\end{lemma}

\begin{proof}
Using \fl(pmn) and \lem{M0} we get:
\begin{equation*}
\begin{split}
\big(\gamma_{-1}\ensquare{-m-1}\gamma_1\big)(n) 
=&\sum_{i\le -m-1}(-1)^m \binom{m+i}m \gamma_{-1}(i)\gamma_1(n-m-1-i)\\
&- \sum_{i\ge0}\,(-1)^m \binom{m+i}m \gamma_1(n-m-1-i)\gamma_{-1}(i)\\
=&\ (-1)^m \sum_{i\in\Z}\binom{m+i}m \hat{e}_{i,i-m+n} = p_m(n)
\end{split}
\end{equation*}
\end{proof}

Now let us apply the construction of \sec{lattice} to the lattice
$\Lambda = \Z$, generated by a single vector $\alpha$, such that
$(\alpha|\alpha)=1$. Then the formula \fl(vanvb) implies that the
elements $v_{\pm\alpha}\in V_\Z$ generate a conformal superalgebra
isomorphic to the Clifford algebra $\goth{Cl}$. As a result we get 
that the vertex algebra $V$ is canonically isomorphic to the vertex 
algebra $V_\Z$ corresponding to the lattice $\Z$ constructed in \sec{lattice}.
This statement is known as bozon-fermion correspondence.

Let us describe the image of the algebra $\^{\goth W}$ in $V_0$ in
terms of the lattice construction. It could be easily proved that 
$p_0 = -\~\alpha \in \goth H \subset V_0$, \ 
$p_1 = \frac12\, \~\alpha\ensquare{-1} \~\alpha - \frac12 D\~\alpha = \omega -
\frac12 D\~\alpha$ (see \sec{lattice}) so that 
$p_1 \ensquare 0 v = Dv$ for all $v\in V$. In general, for $n\ge2$ we have

\begin{align}\label{fl:WinV0}
&p_n = \frac{(-1)^n}{(n+1)!}\times\\[-5pt]
&\begin{aligned}
\scriptstyle{n-2\ \text{times}}
\hskip15pt& \hskip 80pt
\scriptstyle{n\ \text{times}}\\[-10pt]
\Bigg(\binom{n+1}2 
(\!\cdot\hskip-2pt\cdot\hskip-2pt\cdot(\~\alpha\ensquare{-2}\~\alpha)
\overbrace{\ensquare{-1}\~\alpha\ \cdots)\ensquare{-1}\~\alpha}
&\ -\ (\!\cdot\hskip-2pt\cdot\hskip-2pt\cdot(\~\alpha
\overbrace{\ensquare{-1}\~\alpha)
\ensquare{-1}\~\alpha\ \cdots)\ensquare{-1}\~\alpha}
\Bigg).
\end{aligned}\notag
\end{align}

%% file: v0.tex
\subsection{Inside the Heisenberg vertex algebra}\label{sec:v0}
Here we investigate further the embedding $\^{\goth W}\subset V_0$
of the conformal algebra $\^{\goth W}$ into the vertex algebra 
$V_0$, constructed in \sec{fermion}.

\begin{theorem}\label{thm:maximal}\sl
The conformal  algebra $\^{\goth W}\subset V_0$ is a maximal 
conformal subalgebra of $V_0$.
\end{theorem}

By \lem{M0}, the space $V_0$ is a module over the Lie algebra
$\^{gl}_\infty$, and in fact over $\^M$. 
For the proof of \thm{maximal} we need to study the $\^{gl}_\infty$-module
structure of $V_0$. Recall that the Lie algebra $\^{gl}_\infty$ has a
structure quite similar to the structure of a Kac-Moody Lie algebra. 
Let $\goth d \subset \^{gl}_\infty$ be the maximal toral
subalgebra of $\^{gl}_\infty$ spanned by all diagonal matrices and $\c$. 
Let $\goth d' \subset \goth d^*$ be the space of functionals on 
$\goth d$ which take only finitely many nonzero values on $E_{ii}$ 
for $i\in\Z$. Let $\lambda_i\in \goth d', \ i\in\Z$, be the functional
on $\goth d$ such that $\lambda_i(E_{jj})=\delta_{ij}, \
\lambda_i(\c)=0$ and let $\lambda_\c \in \goth d'$ be such that 
$\lambda_\c(E_{jj})=0$ for all $j\in\Z$ and $\lambda_\c(\c)=1$. 
The algebra $\^{gl}_\infty$ is $\goth d$-diagonalizable, the root
vectors being the elements $E_{ij},\ i\neq j,$ whose weights 
$\lambda_{ij}=\lambda_i-\lambda_j\in\goth d'$ are called the roots of
$\^{gl}_\infty$. If $i<j$ we call the root $\lambda_{ij}$ positive, 
otherwise we call it negative. 

The element
\begin{equation}\label{fl:Hij}
\goth d \ni H_{ij}=\ad{E_{ij}}{E_{ji}} = E_{ii}-E_{jj} + 
\begin{cases}
\c&\text{if}\quad j<0 \ \& \  i\ge0\\
-\c&\text{if}\quad j\ge0 \ \&\  i<0\\
0&\text{otherwise}
\end{cases}
\end{equation}
is called a coroot. Denote by 
$\Pi = \set{\lambda\in \goth d'}{\lambda(H_{ij})\in \Z \ 
\text{for all coroots}\ H_{ij}}$ the set of all {\it integral weights}. 

Let $U$ be a  $\^{gl}_\infty$-module. It is  called 
$\goth d$-{\it diagonalizable} if 
$U=\bigoplus_{\lambda\in \goth d'} U_\lambda$, where 
$U_\lambda =$ \break $\set{v\in U}{Hv=\lambda(H)\,v \ \forall\, H\in \goth d}$.
The module $U$ is called {\it integrable} if it is 
$\goth d$-diagonalizable and all $E_{ij}$ for $i\neq j$ are locally nilpotent.
Finally, $U$ is called a {\it lowest weight} module with {\it lowest weight}
$\lambda\in\goth d'$ if it is generated
by a single vector $v\in U_\lambda$ such that $E_{ij}v=0$ when $i>j$
and for any $h\in \goth d$ one has $hv = \lambda(h)v$.  
A lowest weight module $U$ of lowest weight $\lambda$ is integrable if
and only if $\lambda\in\Pi$  and $\lambda(H_{ij})\le0$ when $i<j$,
see \cite[Chapter 10]{kac1}.  For a $\goth d$-diagonalizable module
$U$ denote by $\Xi(U)= \set{\lambda\in\goth d'}{U_\lambda \neq 0}$ 
the set of weights of $U$.

The $\^{gl}_\infty$-module $V_0$ is an integrable irreducible lowest
weight module 
generated by the lowest weight vector $v_0=\1$ of weight $\lambda_c$.
Using the general theory of integrable modules over Kac-Moody
algebras \cite{kac1}, we can easily describe the set of weights 
$\Xi(V_0)$. Before we do that we need a notion from combinatorics,
see \cite{macdonald}.

Let $\kappa = (k_1\ge k_2\ge \ldots) \in \cal P(m)$ be a partition of an
integer $m$ and let $\kappa' = (k'_1\ge k'_2\ge \ldots) \in \cal P(m)$ be
the dual partition, i.e. corresponding to the transposed Young diagram.
Let $l=l(\kappa)$ be the number of rectangles at the main diagonal of the Young
diagrams of $\kappa$ and $\kappa'$. 
Then the pair  of sequences 
$$
\xi = (k_1,k_2-1,k_3-2,\ldots, k_l-l+1),\quad
\eta = (k'_1-1, k'_2-2, k'_3-3,\ldots, k'_l-l)
$$ 
are called Frobenius coordinates of $\kappa$. 
We have 
$$
\xi_1>\xi_2>\ldots>\xi_l>0, \quad 
\eta_1>\eta_2>\ldots>\eta_l\ge0,\quad 
\sum_{i=1}^l \xi_i + \sum_{i=1}^l \eta_i = m.
$$
The Frobenius coordinates $\xi, \ \eta$ of $\kappa$ determine the 
partition $\kappa$ uniquely. We will write $\kappa = \<\xi|\eta\>$.

\begin{lemma}\label{lem:v0}\sl
$\Xi(V_0) = \displaystyle{\bigcup_{m\in\Z_+} \Xi_m}$, where
$$
\Xi_m=\left\{\left.\ \mu(\kappa)=\lambda_\c+
\sum_{j=1}^{l(\kappa)} \big(\lambda_{-\xi_j}-\lambda_{\eta_j}\big)\ 
\right|\ \kappa=\<\xi|\eta\>\in \cal P(m)\right\}.
$$
 The homogeneous component
$V_{0,m}$ of $V_0$ is decomposed into a direct sum of 1-dimensional 
root spaces 
\begin{equation}\label{fl:rsd}
V_{0,m} =\bigoplus_{\lambda \in \Xi_m} V_{0,\lambda},\qquad
\dim V_{0,\lambda}=1.
\end{equation}
\end{lemma}
\begin{remark}
Though this lemma could be proved using only the fact that $V_0$ is
the irreducible lowest weight $\^{gl}_\infty$-module of weight
$\lambda_\c$, in our case the proof is even simpler since we already
know that $\dim V_{0,m} = P(m)$, so we only have to check that weights
$\mu(\kappa)$ indeed appear in $\Xi(V_0)$.
\end{remark}

We will write $l(\mu)$ instead of $l(\kappa)$ if 
$\mu=\mu(\kappa)\in \Xi(V_0)$ is the
weight of $V_0$
corresponding to a partition $\kappa\in\cal P(m)$.

Next we will study the action of the
elementary matrices $E_{ij}\in \^{gl}_\infty$ on $V_0$ in greater
detail. Recall \cite[Chapter 9]{kac1} that there is a contravariant
form $(\,\cdot\,|\,\cdot\,)$ on $V_0$ such that $(\1|\1)=1$, \ 
$(V_{0,m}|V_{0,n})=0$ for $m\neq n$ and $(Au|v) = (u|A^tv)$ for any 
$u,v\in V_0$, $A\in\^M$, $A^t$ being the transposed matrix of
$A$, in particular $\alpha(n)^t = \alpha(-n)$ for 
$\alpha(n)\in H$. The following lemma is proved along the same lines 
as the analogous result about the integrable modules for Kac-Moody
algebras, see \cite[Chapter 10]{kac1}.

\begin{lemma}\label{lem:Eij}\sl
Let $u\in V_{0,\mu}$ be a homogeneous vector of weight
$\mu\in\Xi(V_0)$. If $\op{sign}(i-j)\, \mu(H_{ij})>0$ then
$0\neq E_{ij}u \in V_{0,\mu+\lambda_i-\lambda_j}$ and
$(E_{ij}u|E_{ij}u)=(u|u)$, otherwise $E_{ij}u=0$.
\end{lemma}
Here $H_{ij}$ is given by \fl(Hij). 
Note that from \lem{v0} follows that 
$\mu(E_{ij})\in\left\{0,\pm1\right\}$.

We now investigate the structure of $V_0$ as a $\^{\goth W}$-module.
By definition (see \sec{cm}) this is the same as the action of 
$\^W_+\oplus \Bbbk D$ on $V_0$.
Since $p_1(0)$ acts as $D$, the action of $\^{\goth W}$ on $V_0$
amounts only to the action of the Lie algebra 
$W_+ =  \^W_+ = \Bbbk\<\,p,t\,|\,\ad tp=1\,\>^{(-)}$.

\begin{lemma}\label{lem:W+}\sl
\begin{enumerate}
\item\label{W:hom}
Any $W_+$-submodule of $V_0$ is homogeneous with respect to the root
space decomposition \fl(rsd).
\item\label{W:sum}
Let $v\in V_{0,\lambda}$ for some $\lambda \in \Xi(V_0)$. Then 
$$
W_+v = \bigoplus_{\mu\in\Xi(V_0), \ l(\mu)\le l(\lambda)} V_{0,\mu}.
$$
\end{enumerate}
\end{lemma}

\begin{proof}
Recall that
$V_0$ is a module over the Lie algebra $\^M$, which is a central
extension of the Lie algebra $M$ of infinite matrices with finitely
many non-zero diagonals, and that $\^W$ is embedded in $\^M$ by
formulas \fl(pmn). Let $M_0$ be subalgebra of diagonal matrices of
$M$. Since $\phi(M_0,M_0)=0$, we get that $M_0$ is a subalgebra of
$\^M$. Any functional from $\goth d'$ takes values on $M_0$ as well,
so we have $\goth d' \subset M_0^*$. 
Let $W_0 = \op{Span}\set{p_m(m)}{m\in\Z}$ be the
subalgebra of $W$ consisting of all elements of degree 0. 
We have $W_0\subset M_0$ under the mapping \fl(pmn). 
To prove \ref{W:hom} it is enough to show that any two different weights 
$\lambda, \mu \in \Xi(V_0)$ remain different after restriction to $W_0$.
But this follows from the fact that $\lambda_i$ for $i\in \Z$ are
linearly independent on $M_0$ because 
$\lambda_i\big(p_m(m)\big) = (-1)^m \binom{i+m}m$.

For the proof of \ref{W:sum} we note that by \fl(pmn) every element
$p^mt^n\in W_+$ is an infinite linear combination of $E_{ij}$'s
such that either $i\ge0$ or $j<0$. Therefore, by \lem{Eij} and 
\ref{W:hom} we get that if 
$\mu \in \Xi(V_0)$ then for any $i,j\in\Z$ such that either $i\ge0$ or 
$j<0$ the weight $\mu+\lambda_i-\lambda_j$ appears in $W_+V_\mu$. 
Every weight $\lambda$ such that $l(\lambda)\le l(\mu)$ could be
obtained by a successive application of this operation but no weight
of length more than $l(\mu)$ can be obtained.
\end{proof}

In particular all weights of $\^{\goth W}\subset V_0$ are of length 1.

\begin{proof}[Proof of \thm{maximal}]
Let $v\in V_0\ssm \^{\goth W}$. We have to prove that the
conformal algebra generated by $v$ and $\^{\goth W}$ is the whole
$V_0$. By \lem{W+} we can assume that $v$ is homogeneous of some
weight $\mu \in \Xi(V_0)$ such that $l(\mu)>1$. The only
weight of degree 4 which is not in $\Xi(\^{\goth W})$ is 
$\lambda_{-2}+\lambda_{-1}-\lambda_0-\lambda_1$ which is of length 2,
hence by \lem{W+}\ref{W:sum},  
$V_{0,4}\subset W_+ v + \^{\goth W}$. Therefore it is enough to
prove that $V_0$ is generated as a conformal algebra by 
$\^{\goth W}$ and any element $u\in V_{0,4}\ssm \^{\goth W}$. Take 
$u=\~\alpha \ensquare{-3} \~\alpha$. 

Let $\goth L\subset V_0$ be the conformal algebra generated by 
$\^{\goth W}$ and $u$. Assume on the contrary that $\goth L\subsetneq
V_0$. Then \lem{W+}\ref{W:sum} implies that there is an integer $l\ge 3$ such
that $\Xi(\goth L)$ does not contain weights of length $l$. Let $l$ be
the minimal possible among such integers.  Let 
$$
\mu = \lambda_{-l}+\ldots+\lambda_{-3}+\lambda_{-2}
-\lambda_2-\lambda_3-\ldots-\lambda_l+\lambda_{\c} \in \Xi(\goth L),
$$
$l(\mu)=l-1$, 
and let $v\in \goth L_{\mu}\subset V_{0,\mu}$. Let 
$\mu' = \mu+\lambda_{-1}-\lambda_0$ so that $l(\mu')=l$ and hence
$\mu'\not\in \Xi(\goth L)$. We prove that the projection of $u\ensquare 2 v$
onto  $V_{\mu'}$ is non-zero, therefore $\goth L_{\mu'}\neq0$ so we
arrive to a contradiction.

Using that $\~\alpha(k) = -t^k = -\sum_{i\in \Z} E_{i,i-k}$
and \fl(verjac), we get 
\begin{equation*}
\big(\~\alpha\ensquare{-3}\~\alpha\big)\ensquare 2 v 
=\sum_{s\ge0}\(\binom s2+\binom {s+2}2\)
\sum_{i,j\in\Z}E_{i,i+s+1}E_{j,j-s}\,v. 
\end{equation*}
The only pairs of $E_{ij}$'s in the above expression 
that do not kill $v$ and move $v$ to
$V_{\mu'}$ are
$E_{-1,i}E_{i,0}$ for $2\le i\le l$,\ 
$E_{i,0}E_{-1,i}$ for $-l\le i\le -2$, and
$E_{-1,0}E_{i,i}$ for either $2\le i\le l$ or $-l\le i\le -2$.
Let $w = E_{-1,2}E_{2,0}\,v$. 
Using \lem{Eij} we get that 
$E_{-1,i}E_{i,0}\,v=w$ for all $2\le i\le l$ and  
$E_{i,0}E_{-1,i}\,v=-w$ for all $-l\le i\le -2$. Also, 
$E_{-1,0}E_{i,i}\,v = - E_{-1,0}E_{-i,-i}\,v$. Summing up, we get 
\begin{align*}
\op{Pr}_{\,V_{\mu'}}\big(\~\alpha\ensquare{-3}\~\alpha\big)\ensquare 2 v 
& = \sum_{s=2}^l\(\binom s2+\binom{s+2}2\) w  \ - \ 
\sum_{s=1}^{l-1}\(\binom s2+\binom{s+2}2\) w\\[3pt]
&= \(\binom l2 + \binom {l+2}2 - 3\)x \neq0
\end{align*}
\end{proof}

\subsection{Representation theory of $\^{\goth W}$}
\label{sec:repW}
The technique of \sec{v0} allows to investigate the 
$\^{\goth W}$-module structure of other spaces 
$V_\lambda$. This section deviates from the main exposition and will
not be used later in this paper. We assume here  that $\Bbbk$ is
an algebraically closed field of characteristic 0.

Let $\Lambda=\Z\beta$ be a lattice of rank 1 generated by a single
vector $\beta$. As before let $\goth h=\Lambda\otimes \Bbbk$.  Let
$V=V_\Lambda=\bigoplus_{i\in\Z}V_{i\beta}$ be the corresponding vertex
algebra constructed in \sec{lattice}. Take $\alpha =
\frac\beta{\sqrt{(\beta|\beta)}}\in\goth h$.  Then $(\alpha|\alpha)=1$
and the field $\~\alpha=\alpha(-1)\1\in V_0$ generates a copy of the
Heisenberg conformal algebra $\goth H$ and also gives an embedding of
$\^{\goth W}$ in $V_0$ by formulas
\fl(WinV0). So $V_{i\beta}$ becomes a module over  $\^{\goth W}$ and
over $W_+$.

If $(\beta|\beta)\neq1$, then it is not difficult to see that
$V_{i\beta}$ is an irreducible $W_+$-module. In fact, in this case we
have $W_+v_{i\beta}=V_{i\beta}$. If $(\beta|\beta)=1$, and in this
case $\Lambda = \Z$, then as in \sec{v0} we have an action of 
$\^{gl}_\infty$ on $V_{i\beta}=V_i$. The module $V_i$ is an
irreducible integrable lowest weight $\^{gl}_\infty$-module of the
lowest weight $\lambda_\c - \lambda_0-\lambda_1-\ldots-\lambda_{i-1}$
if $i\ge0$ and $\lambda_\c + \lambda_{-1}+\ldots+\lambda_i$ if $i<0$.

We need to introduce the biased Frobenius coordinates of a partition. 
Let $\kappa = (k_1\ge k_2\ge \ldots) \in \cal P(m)$ be a partition of an
integer $m$ and let $\kappa' = (k'_1\ge k'_2\ge \ldots) \in \cal P(m)$ be
the dual partition. Denote by $l_i(\kappa)$ the number of squares on
the $i$th diagonal of the Young diagram of $\kappa$, so that 
$l_i(\kappa)=l_{-i}(\kappa')$. The biased Frobenius coordinates of 
$\kappa$ are sequences
$$
\xi = (k_1-i,k_2-i-1,\ldots,k_{l_i}-i-l_i+1),\quad
\eta = (k'_1+i-1, k'_2+i-2,\ldots,k'_{l_i+i}-l_i).
$$
They determine $\kappa$ uniquely and we will write 
$\kappa = \<\xi|\eta\>_i$.

Let $\Xi(V_i)\subset \goth d'$ be the set of weights of the module
$V_i$. In analogy with \lem{v0} we have that   
$\Xi(V_i) = \displaystyle{\bigcup_{m\in\Z_+} \Xi_m(V_i)}$, where
$$
\Xi_m(V_i)=\left\{\left.\ \mu_i(\kappa)=\lambda_\c+
\sum_{j=1}^{l_i(\kappa)} \lambda_{-\xi_j}-
\sum_{j=1}^{l_i(\kappa)+i}\lambda_{\eta_j}\ 
\right|\ \kappa=\<\xi|\eta\>_i\in \cal P(m)\right\}
$$
is the set of weights appearing in the homogeneous component 
of $V_i$ consisting of elements of degree $m+\frac12i^2$. As before we
have $\dim V_{i,\mu}=1$ for every $\mu\in\Xi(V_i)$.

Both statements \ref{W:hom} and \ref{W:sum} of \lem{W+} remain without changes.
In particular, dimensions of homogeneous components of
$W_+$-submodules of $V_i$ grow polynomially. 

\begin{question}
Is it true that all lowest weight 
$W_+$-modules of polynomial growth are obtained in this way?
\end{question}

%% file: rootsys.tex
\subsection{Conformal subalgebras of lattice vertex superalgebras}
\label{sec:rootsys}
Let $\Lambda$ be an integer lattice and let 
$V_\Lambda= \bigoplus_{\lambda\in\Lambda}V_\lambda$ be the lattice 
vertex superalgebra constructed in \sec{lattice}. Recall that $V_0$
contains the Heisenberg conformal algebra $\goth H$, which is 
spanned over $\Bbbk[D]$ by elements of the form $\~h$ for 
$h\in \goth h = \Bbbk\otimes \Lambda$ and $\1$. 
The left regular action of this
elements is given by $Y\big(\~h\big) = \sum_{n\in\Z} h(n)\, z^{-n-1}$, where
the operators $h(n)$ define a representation of the Heisenberg algebra
$H$ on $V_\Lambda$. 

Let $\goth L\subset V_\Lambda$ be a conformal subalgebra of
$V_\Lambda$. Assume that $\goth L$ is homogeneous with respect to the
grading by $\Lambda$, that is, $\goth L =
\bigoplus_{\lambda\in\Lambda}\goth L_\lambda$, where $\goth L_\lambda
= \goth L \cap V_\lambda$. The set $\Delta = \set{\lambda\in
\Lambda\ssm 0}{\goth L_\lambda \neq 0}$ is called the {\it root
system} of $\goth L$. We will always assume that $\Delta =
-\Delta$. This happens, for example, if $\goth L$ is closed under the
involution $\sigma:V_\Lambda \to V_\Lambda$, corresponding to the
automorphism $\lambda \mapsto -\lambda$ of the lattice $\Lambda$. 

We will also assume that $\goth L$ is closed under the action of the
conformal algebra $\goth H$. In this case it is easy to show that
$\goth L$ must contain the elements $v_\lambda$ for each $\lambda \in
\Delta$. Therefore, the Lie conformal superalgebra $\goth L'\subset
V_\Lambda$ generated by the set $\set{v_\lambda}{\lambda\in \Delta}$
will be a subalgebra of $\goth L$ and will be graded by the same root
system.  
From the formula \fl(vanvb) for the products in $V_\Lambda$ follows
that if $\alpha, \beta \in \Delta$ be such that $k=(\alpha|\beta)<0$
then $v_\alpha\ensquare{-k-1}v_\beta = \pm v_{\alpha+\beta}\in \goth
L'$, hence $\alpha+\beta \in \Delta$. 

The following proposition summarizes properties of root systems.

\begin{proposition}\label{prop:rs}\sl
\begin{enumerate}
\item
A set $\Delta \subset \Lambda$ is a root system if and only if the
Lie conformal superalgebra generated by the set 
$\set{v_\lambda}{\lambda\in \Delta}$ does not contain any homogeneous
components other than $\goth L_\lambda$ for $\lambda \in \Delta$ and 
$\goth L_0$.
\item
If $\Delta\subset \Lambda$ is a root system, then $\Delta$ is closed
under the negation $\lambda \mapsto -\lambda$ and under the partial
summation: 
\begin{equation}\label{fl:partsum}
\alpha, \beta\in \Delta, \ \ (\alpha|\beta)<0 \ \ \Longrightarrow \ \ 
\alpha + \beta \in \Delta
\end{equation}
\end{enumerate}
\end{proposition}

We are mostly interested in the case when the root system $\Delta$ is
finite. If $\Delta$ contains a vector $\lambda$ such that
$(\lambda|\lambda)<0$ then by \prop{rs}(b) $k \lambda \in \Delta$ for
all $k = 1,2,\ldots$. The following lemma suggests that we should
restrict ourselves to the case when the form $(\,\cdot\,|\,\cdot\,)$
is semi-positive definite.

\begin{lemma}\label{lem:rs}\sl
Let $\Delta \in \Lambda$ be a set 
 closed under the partial summation \fl(partsum) such that 
$\Delta =-\Delta$ and $\op{Span}_\Z \Delta = \Lambda$. 
\begin{enumerate}
\item\label{rs:neg}
Assume that the form $(\,\cdot\,|\,\cdot\,)$ is not semi-positive
definite, i.e., there is $\alpha \in \Lambda$ such that
$(\alpha|\alpha)<0$. Then there is some $\delta \in \Delta$ such that 
$(\delta|\delta)<0$. 
\item\label{rs:spd}
Assume that  the form $(\,\cdot\,|\,\cdot\,)$ is semi-positive but not
positive definite. Then there exists some $\delta \in \Delta$ such that 
$(\delta|\delta)=0$.
\end{enumerate}
\end{lemma}

\begin{proof}
Both statements are proved by a standard argument. Let us prove (\ref{rs:neg}),
the prove of (\ref{rs:spd}) being identical. Assume on the contrary that 
$(\lambda|\lambda)\ge 0$ for every $\lambda \in \Delta$. Let $\alpha
\in \Lambda$ be such  that $(\alpha|\alpha)<0$. Then $\alpha$ could
be expressed as a linear combination $\alpha = k_1 \alpha_1 +\ldots k_l
\alpha_l$ of elements $\alpha_i \in \Delta$ with integer
coefficients. Since $\Delta = -\Delta$ we can assume that all $k_i>0$.
Assume that the above is a combination of this kind with the minimal
possible value of $\sum_i k_i$. Then since
$(\alpha_i|\alpha_i)\ge 0$ for all $i$ and $(\alpha|\alpha)<0$ we must
have $(\alpha_i|\alpha_j)<0$ for some $i\neq j$, but then $\alpha_i +
\alpha_j \in \Delta$ and we could make the expression for $\alpha$
with a smaller sum of coefficients.
\end{proof}

The following sections will be dedicated to the classification of
finite root systems. In case when the bilinear form has a non-trivial
kernel, we allow for a bigger class of root systems. Let 
$\Lambda_0 =\set{\lambda\in\Lambda}{(\lambda|\Lambda)=0}
\subset \Lambda$ be the
sublattice of isotropic vectors. Let $\ol\Lambda = \Lambda/\Lambda_0$
be the positive definite quotient. We will denote the projection of an
element $\lambda \in \Lambda$ to $\ol{\Lambda}$ by $\ol\lambda$.
A set $\Sigma \in \Lambda$ is called 
{\it almost finite} if $\ol\Sigma \subset \ol\Lambda$ is finite.

We will adopt the following notations. Let $\Lambda$ be a
semi-positive definite integer lattice and $\Lambda_0$ and
$\ol\Lambda$ be as above. Let $\Delta \subset \Lambda$ be a root
system. Denote by $\Delta_0 = \Delta \cap \Lambda_0$ the set of
isotropic roots and by $\Delta^\times = \Delta \ssm \Delta_0$ the set
of all real roots. Let $\ol\Delta = (\Delta/\Lambda_0) \ssm \{0\}$ be
the positive definite root system in $\ol\Lambda$ obtained from the
projection of $\Delta$ to $\ol\Lambda$.

We also need the following definition. A root system $\Delta \subset
\Lambda$ is called {\it decomposable} if it can be represented as a
disjoint union $\Delta = \Delta_1 \sqcup \Delta_2$ such that for any 
$\alpha \in \Delta_1, \ \beta \in \Delta_2$ we have $\alpha + \beta
\not\in \Delta$. Otherwise $\Delta$ is called {\it indecomposable}. By
\prop{rs}(b) if 
$\Delta= \Delta_1 \sqcup \Delta_2$ is decomposable then
$(\Delta_1|\Delta_2)=0$ and  the Lie conformal superalgebra  $\goth L$ 
splits into a direct sum
$\goth L = \goth L_1 \oplus \goth L_2$ such that 
$\ad{\goth L_1}{\goth L_2}=0$. In the other direction, however, if
$\goth L= \goth L_1 \oplus \goth L_2$ so that 
$\ad{\goth L_1}{\goth L_2}=0$, and we set $\Delta_i$ to be the
root system of $\goth L_i$, then though $(\Delta_1|\Delta_2)=0$, in general 
$(\Delta_1 + \Delta_2)\cap \Delta \neq \emptyset$ unless the form is
positive definite.

%% file: rank1.tex
\subsection{Rank 1 case}
\label{sec:rank1}
Let $\alpha \in\Lambda$ be a vector in an integer lattice $\Lambda$.
Let $\goth L\subset V_\Lambda$ be the conformal superalgebra
generated by $v_\alpha$ and $v_{-\alpha}$ in the lattice vertex
superalgebra $V_\Lambda$, constructed in \sec{lattice}. The algebra
$\goth L = \bigoplus_{\lambda\in\Lambda}\goth L_\lambda$ is graded by
the lattice $\Lambda$. We start with determining all cases when $\goth
L$ does not contain any homogeneous components other than $\goth
L_{-\alpha}\ni v_{-\alpha}, \ \goth L_0 \ni \1$ and $\goth L_\alpha
\ni v_\alpha$.

Clearly, if $(\alpha|\alpha)<0$, then take $n=-(\alpha|\alpha)-1\ge0$ and
get $v_\alpha \ensquare n v_{-\alpha} = v_{2\alpha}\in \goth L$. hence
all $v_{j\alpha}\in\goth L$ for $j\in\Z$. Also, if 
$(\alpha|\alpha)=0$ then all products in $\goth L$ are 0 so this case
is not interesting. Therefore without a loss of generality we assume
that  $(\alpha|\alpha)>0$.

\begin{proposition}\label{prop:rank1}\sl
If $(\alpha|\alpha)=1$ then $\goth L = \goth{Cl}$ is the Clifford
conformal superalgebra. 

If  $(\alpha|\alpha)=2$
then $\goth L=\goth{Aff}(sl_2)$ is the affine conformal algebra 
$\^{sl}_2$.

If $(\alpha|\alpha)=3$ then $\goth L$ is the central extension
of the $N=2$ simple conformal
superalgebra.

If $(\alpha|\alpha)=4$ then $\goth L$ is the conformal
algebra $\^{\goth K}$, constructed in \sec{weyl}. 

Finally, if  $(\alpha|\alpha)\ge 5$ then
$\goth L = V_{\Z\alpha}$.
\end{proposition}

\begin{proof}
First we show that if $(\alpha|\alpha)=1,2,3$ or $4$ the conformal
algebra $\goth L \subset V_\Lambda$ generated by $v_\alpha$ and
$v_{-\alpha}$ is as claimed. 
 
As it was remarked at the end of \sec{lattice}, all calculations in
vertex algebra $V_\Lambda$ are very explicit. So we just have to read
off the defining relations of conformal superalgebras from the
formula \fl(vanvb) for the products in $V_\Lambda$. Of course, the
case when $(\alpha|\alpha)$ is 1 or 2 is well-known. Let us do the
most difficult case when $(\alpha|\alpha)=4$. In this case we can
identify $\Z\alpha$ with $2\Z\subset \Z$ by letting $\alpha=2$.

By \lem{WinV} the elements $p_m = v_{-1} \ensquare{-m-1} v_1\in V_0$ and $\1$
span a copy of $\^{\goth W}$ over $\Bbbk[D]$ so that the products are
given by \fl(confpp) and \fl(confphiW). Set $u_m = v_{-1} \ensquare{-m-1}
v_{-1}\in V_{-2}, m\ge 1$. Recall then there is an involution
$\sigma:V_\Lambda \to V_\Lambda$ induced by the involution
$\lambda\mapsto -\lambda$ of the lattice $\Lambda$, so that we have 
$\sigma(u_m) = v_1 \ensquare{-m-1} v_1\in V_2$.
We have to show
that the formulas \fl(udepend)--\fl(Phieven) hold for $p_m, u_m$,
$\sigma(u_m) \in V_\Z$ and $\c=\1$ and also $u_m \ensquare k u_n = 
\sigma(u_m)\ensquare k \sigma(u_n) = 0$ for all integer 
$m,n\ge 1, \ k \ge 0$. 

Let us for example check \fl(confpsu). First we note that 
\begin{equation*}
\begin{split}
v_{-1}\ensquare i \bigl(v_1\ensquare{-n-1}v_1\bigr) 
&= -v_1\ensquare{-n-1}\bigl(v_{-1}\ensquare i v_1\bigr) + 
\ad{v_{-1}(s)}{v_1(-n-1)} v_1 \\
&= -\delta_{i,0}\, v_1(-n-1)\1 + \1(i-n-1)\,v_1 \\
&=-\delta_{i,0} \,\frac 1{n!} \,D^n v_1 + \delta_{i,n}\, v_1.
\end{split}
\end{equation*}
Using this we calculate
\begin{equation*}
\begin{split}
p_m \ensquare k \sigma(u_n) &= \bigl(v_{-1}\ensquare{-m-1} v_1\bigr)
\ensquare k \bigl(v_1\ensquare{-n-1}v_1\bigr)\\
&= \sum_{s\le -m-1}\binom{-s-1}m\, v_{-1}(s)\,v_1(k-m-1-s)\,v_1(-n-1)\,v_1  \\
&-(-1)^m 
\sum_{s\ge0} \binom{m+s}m \,v_{-1}(k-m-1-s)\,v_{-1}(s)\,v_1(-n-1)\,v_1.  
\end{split}
\end{equation*}
The first sum here is 0 because $k-m-1-s\ge 0$, hence 
$v_1(k-m-1-s)$ commutes with $v_1(-n-1)$ and $v_1(k-m-1-s)\, v_1 = v_1
\ensquare{k-m-1-s}v_1 =0$. The second sum gives
\begin{equation*}
\frac{(-1)^m}{n!}\, v_1\ensquare{k-m-1} \big(D^n v_1\big) - 
(-1)^m \binom{m+n}m\, v_1\ensquare{k-m-1-n}v_1,
\end{equation*}
and \fl(confpsu) follows. The other formulas are checked in the same
way.  Instead of checking \fl(udepend) directly, we note that by
\sec{repW} the $\^{\goth W}$-modules  
$U(W_+)\, v_{\pm 2}$ and $\goth J$ are both  irreducible highest
weight $W_+$-modules
corresponding to the same highest weight, hence they must be
isomorphic.  

The verification of the conformal cocycle formula \fl(Phieven) is
done in the same way.

We are left to show that if $n = (\alpha|\alpha)\ge 5$, then $\goth L =
V_\Lambda$. Recall from \sec{fermion} that we have an embedding
$\^{\goth W}\subset V_0$ given by \fl(WinV0) such that 
$\^{\goth W}_i = \op{Span}_\Bbbk
\left\{p_{i-1},Dp_{i-2},\ldots,D^{i-1}p_0\right\}\subset V_{0i}$ 
for $i\ge1$ so that $\^{\goth W}_i =V_{0i}$ for $0\le i\le 3$.
Simple calculations
show that $v_{\alpha}\ensquare{n-i}v_\alpha \in \^{\goth W}_{i-1}$
for $1\le i \le4$ and also 
$\left\{\1,p_0,p_1,p_2\right\}\subset 
\op{Span}_\Bbbk\set{v_{-\alpha}\ensquare{n-i}v_\alpha}{1\le i\le4}$.
Since $\^{\goth W}$ is generated by $\left\{\1,p_0,p_1,p_2\right\}$, 
we have that $\^{\goth W}\subset \goth L$. However, 
$v_{\alpha}\ensquare{n-5}v_\alpha \in V_{04}\ssm \^{\goth W}_4$,
therefore, since $\^{\goth W}$ is a maximal conformal subalgebra in
$V_0$ by \thm{maximal}, we must have $\goth L = V_0$. 
\end{proof}

It follows that all possible finite root systems $\Delta$ of rank 1
are from the following list:

\smallskip
\begin{tabular}{rl}
$A_1$:&
$\Delta= \{\pm\alpha\}, \ (\alpha|\alpha)=2$;\\
$B_1$:&
$\Delta=\{\pm\alpha\}, \ (\alpha|\alpha)=1$;\\
$C_1$:&
$\Delta=\{\pm\alpha\}, \ (\alpha|\alpha)=4$;\\
\end{tabular}

\begin{tabular}{rl}
$BC_1$:&
$\Delta=\{\pm\alpha, \, \pm 2\alpha\}, \ (\alpha|\alpha)=1$;\\
$B'_1$:&
$\Delta=\{\pm\alpha\}, \ (\alpha|\alpha)=3$.
\end{tabular}
 
\smallskip
We will generalize this result for the case of a higher ranking
integer positive definite lattice $\Lambda$ and finite root system
$\Delta\in \Lambda$ in \sec{posdef}, see \thm{pd}.   

Of special interest is the case when we take the root system 
$\Delta = \{\,\pm1, \pm2\,\}\subset \Z$ is of type $BC_1$. 
Then the  conformal superalgebra 
$\goth L\subset V_\Z$ generated by the set $\{\,v_{\pm1}, v_{\pm2}\,\}$ 
is isomorphic to an extension of $\goth K$ by the Clifford 
conformal superalgebra $\goth{Cl}$ 
such that $\goth L\even = \goth L_{-2}
\oplus \goth L_0 \oplus \goth L_2 =\^{\goth K}$ and 
$\goth L_{-1}\oplus \goth L_1 \oplus \1 = \goth{Cl}$. 

We also remark that  $\goth L$ and $\goth K$ are maximal among 
conformal subalgebras of $V_\Z$ graded by the corresponding root system.

%% file: rank2.tex
\subsection{Rank 2 case}
\label{sec:rank2}
Consider now two vectors $\alpha,\beta\in \Lambda$, where $\Lambda$
is an integer lattice as before. Let $V_\Lambda$ be the vertex
superalgebra corresponding to $\Lambda$, and let
$\goth L\subset V_\Lambda$
be the conformal  superalgebra generated by $v_{\pm\alpha}$ and 
$v_{\pm\beta}$. Since the generators of $\goth L$ are
homogeneous, $\goth L = \bigoplus_{\lambda\in\Lambda} \goth L_\lambda$
is graded by $\Lambda$. Let 
$\Delta = \set{\lambda\in\Lambda\ssm\{0\}}{\goth L_\lambda\neq0}$ be the root
system of $\goth L$. If $(\alpha|\beta)=0$ then 
$\goth L = \goth L_1 \oplus \goth L_2$ is decomposed into a direct sum
of ideals $\goth L_1 = \<v_{\pm\alpha}\>$ and $\goth L_2=\<v_{\pm\beta}\>$ and 
$\Delta = \Delta_1 \sqcup \Delta_2$ where $\Delta_1 \subset \Z\alpha$
and $\Delta_2 \subset \Z\beta$, so everything is reduced to the case
of \sec{rank1}. Therefore without a loss of generality we assume that 
$(\alpha|\beta)<0$.

Now we formulate an analogue of \prop{rank1} for the case of two
vectors.

\begin{proposition}\label{prop:rank2}\sl
Let $\Lambda = \Z\alpha + \Z\beta$ be an integer lattice of rank 2. 
Assume that the  conformal superalgebra $\goth L\subset V_\Lambda$ 
generated by $v_{\pm\alpha}$ and $v_{\pm\beta}$ is graded by a finite or an
almost finite root system $\Delta \subset \Lambda$.
Then there are only the following possibilities:

\smallskip
\begin{tabular}{rl}
\rm (i)& 
$(\alpha|\alpha)=(\beta|\beta)=2, \ \ (\alpha|\beta)=-1$ \\
\rm (ii)& 
$(\alpha|\alpha)=2, \ \ (\beta|\beta)=1, \ (\alpha|\beta)=-1$\\
\rm (iii)& 
$(\alpha|\alpha)=4, \ \ (\beta|\beta)=2, \ (\alpha|\beta)=-2$\\
\rm (iv)& 
$(\alpha|\alpha)=(\beta|\beta)=1, \ \ (\alpha|\beta)=-1$ \\
\rm (v)& 
$(\alpha|\alpha)=(\beta|\beta)=2, \ \ (\alpha|\beta)=-2$ \\
\rm (vi)& 
$(\alpha|\alpha)=(\beta|\beta)=3, \ \ (\alpha|\beta)=-3$ \\
\rm (vii)& 
$(\alpha|\alpha)=(\beta|\beta)=4, \ \ (\alpha|\beta)=-4$ \\
\rm (viii)& 
$(\alpha|\alpha)=4, \ \ (\beta|\beta)=1, \ (\alpha|\beta)=-2$
\end{tabular}
\smallskip

In cases {\rm(i)--(iii)} $\Lambda$ is positive definite, 
in {\rm (iv)--(viii)} $\Lambda$ is semi-positive definite with the
kernel of the bilinear form spanned by single vector $\delta$ 
given by $\delta=\alpha+\beta$ in cases {\rm(iv)--(vii)} and  
$\delta=\alpha+2\beta$ in case {\rm(viii)}. The root system is 
$$
\Delta = \begin{cases}
\{\,\pm\alpha,\ \pm\beta,\ \pm(\alpha+\beta)\,\}&\text{in cases {\rm(i),
(ii)} and \rm(iv)}\\
\{\,\pm\alpha,\ \pm\beta,\ \pm(\alpha+\beta),\
\pm(\alpha+2\beta)\,\}&\text{in case \rm (iii)}\\
\set{k\delta,\pm\alpha+k\delta,\ \pm\beta+k\delta}{k\in\Z}&\text{in
cases \rm (v)--(viii)}\\
\end{cases}
$$
If $(\lambda|\lambda)=1$ or 2 for $\lambda\in \Delta$ then 
$\goth L_\lambda \cong \Bbbk[D]v_\lambda$, if $(\lambda|\lambda)=3$ or 4 
then $\op{rk}_{\Bbbk[D]}\goth L_\lambda = \infty$. For an isotropic
$\lambda \in \Delta$ we have $\op{rk}_{\Bbbk[D]}\goth L_\lambda =
1,2,3,\infty,\infty$ in the cases {\rm (iv)--(viii)} respectively. 
Finally, $\op{rk}_{\Bbbk[D]}\goth L_0 =
2,1,\infty,0,2,4,\infty,\infty$ in the cases {\rm (i)--(viii)} respectively. 
\end{proposition}
\begin{proof}
By \prop{rank1} the square lengths of $\alpha$ and $\beta$ could be
only 1, 2, 3 or 4. Also, 
$(\alpha|\beta)\ge - \frac12 \big((\alpha|\alpha)+(\beta|\beta)\big)$,
because otherwise we would have $(\alpha+\beta|\alpha+\beta)<0$. So we
have only finitely many possibilities. One can easily check using
formulas \fl(vanvb) that every choice of vectors $\alpha$ and
$\beta$, not listed in (i)--(viii) will give an infinite root system. 
So it remains to show that in each of the cases (i)-(viii) the conformal
algebra $\goth L$ generated by $\{\,v_{\pm\alpha},\, v_{\pm\beta}\,\}$
will be in fact graded by the
corresponding  root system $\Delta$. Let us
check this for the most difficult case (iii). In this case 
$\goth L_0 = \^{\goth W}_1\oplus \^{\goth W}_2$ is a direct sum of two 
copies of $\^{\goth W}$, corresponding to vectors $\frac12\alpha$ and 
$\frac12\alpha+\beta$, so that the subalgebras $\goth L_\alpha \oplus
\^{\goth W}_1 \oplus \goth L_{-\alpha}$ and $\goth L_{\alpha +2\beta}\oplus
\^{\goth W}_2 \oplus \goth L_{-\alpha-2\beta}$ of $\goth L$ are
isomorphic to the conformal algebra $\^{\goth K}$ constructed in \S2. So the
claim follows from \prop{rank1}.
\end{proof}

\subsection{Case when the bilinear form is positive definite}
\label{sec:posdef}
\begin{theorem}\label{thm:pd}\sl
Assume $\Lambda$ is a positive definite integer lattice and 
$\Delta \subset \Lambda$ is a finite indecomposable root system of some 
 conformal superalgebra 
$\goth L = \bigoplus_{\lambda\in\Delta}\goth L_\lambda
\subset V_\Lambda$, such that $\Delta =-\Delta$ and $\goth L$ is
generated by the set $\set{v_\lambda}{\lambda \in 
\Delta}$. Then there are only the following possibilities:

\noindent
\begin{flushright}
\begin{tabular}{rp{3.9in}}
$A$-$D$-$E$:& $\Delta$ is a simply-laced finite Cartan root
system of type 
$A_n$, $n\ge 1$, $D_n$, $n\ge 4$ or $E_n$, $n=6,7,8$ such that
$(\lambda|\lambda)=2$ for all roots $\lambda \in \Delta$.\\
$B$:&
$\Delta$ is a finite Cartan root system of type
$B_n, n\ge 1$, the short roots have square length 1 and long roots
have square length 2. (When $n=1$, \ $\Delta = \{\pm1\} \subset \Lambda=\Z$). \\
$C$:&
$\Delta$ is a finite Cartan root system of type
$C_n, n\ge 1$, the short roots have square length 2 and long roots
have square length 4. (When $n=1$, \ $\Delta = \{\pm2\} \subset \Lambda=\Z$). \\
$BC$:&
$\Delta$ is the union of $B_n$ and $C_n$ for $n\ge 1$.\\
\end{tabular}

\noindent
\begin{tabular}{rp{3.9in}}
$B^0$:&
$\Delta$ is a subset of $B_n$ consisting of all the short roots of
$B_n$ and half of the long roots: if $\alpha_1, \ldots, \alpha_n$ is
the basis of $\Lambda$ consisting of short roots so that
$(\alpha_1|\alpha_j)=0$ then all the long roots of $\Delta$ are of the
form $\alpha_i - \alpha_j, \ i\neq j$.\\
$B'_1$:&
$\Delta = \{\pm3\}\subset \Lambda=\Z$.
\end{tabular}
\end{flushright}
\end{theorem}

We will call a vector of
square length 1 short, of square length 2 long and of square length 4
extra-long. 

\begin{proof}
Let $\alpha,\beta \in \Delta$ be a pair of roots. The root system
which they generate must be of one of the types (i)--(iii) from
\prop{rank2}. It follows that the Cartan number $\<\alpha, \beta\> =
\frac{2 (\alpha|\beta)}{(\alpha|\alpha)}$ is an integer, hence $\Delta$
must lay inside a Cartan root system $\Phi$. Moreover it follows from
the structure of root systems of rank 2 in cases (i)--(iii) that
$\Phi$ cannot be of types $F_4$ and $G_2$, and the condition for
the square lengths of the root vectors hold. 

Next we want to prove that if $\Delta$ is a finite Cartan root system
of the type other than $G_2$ and $F_4$ and the length of roots are as
prescribed by the theorem,
then $\Delta$ is indeed the root system of the 
conformal superalgebra $\goth L\subset V_\Lambda$ generated by the set 
$\set{v_\lambda}{\lambda\in\Delta}$. If $\Delta$ is a
simply-laced root system of types $A$-$D$-$E$ then $\goth L$ is the
affine Kac-Moody conformal algebra, as it is well-known. If 
$\Delta$ is of type $B$ then it is equally easy to check that 
the space 
$\goth L = \bigoplus_{\alpha\in \Delta} \Bbbk[D] v_\alpha\oplus \goth
H \subset V_\Lambda$ will be closed under the products \fl(vanvb),
hence it is the desired subalgebra. 

Let  $\Delta$ be of type $C_n$. Let $\alpha_1, \ldots, \alpha_n\in
\Delta$ be the basis of $\Lambda$ consisting of pairwise orthogonal
extra-long roots. Take 
$\goth L_0 = \^{\goth W}_1 \oplus \ldots \oplus \^{\goth W}_n \subset V_0$, 
where $\^{\goth W}_i$ is the Weyl conformal algebra spanned by 
$v_{\alpha_i}\ensquare n v_{-\alpha_i},\ n\in \Z$, see \lem{WinV}. 
Take $\goth L_{\alpha_i}$ to be equal to the conformal Jordan triple
system, constructed in \S2, so that the subalgebra 
$\goth L_{-\alpha_i}\oplus \^{\goth W}_i \oplus \goth L_{\alpha_i}$ is
isomorphic to the  conformal algebra $\goth K$. 
If $\beta\in \Delta$ is a short root then take $\goth L_\beta =
\Bbbk[D]v_\beta$. Using calculations of \prop{rank2} it is easy to see
that $\goth L = \goth L_0 \oplus\bigoplus_{\alpha\in \Delta}\goth L_\alpha$ 
is closed under the products \fl(vanvb).

Finally, if $\Delta$ is of type $BC$ then the corresponding 
conformal superalgebra $\goth L$ is easily obtained by combining the 
conformal superalgebras corresponding to the subsystems of $\Delta$ of
types $B$ and $C$.

So let $\Delta$ be a finite root system, and let $\Phi\supset \Delta$
be a minimal Cartan root system containing $\Delta$. We prove that if
$\Phi$ is either simply-laced or is of type $C$ or $BC$ then 
$\Delta = \Phi$. 

Assume first that $\Phi$ is simply-laced. Let 
$\alpha \in \Phi\ssm\Delta$. Since $\Delta$
spans $\Lambda$ over $\Z$ we can write $\alpha$ as a linear
combination of elements of $\Delta$ with integer coefficients. Let 
$\alpha = \sum_i k_i \alpha_i,
\ k_i \in \Z, \alpha_i \in \Delta$, be such a linear combination of
the minimal length.  
Since $\Delta = -\Delta$ we can assume that all $k_i>0$.  But then
since $(\alpha|\alpha) = (\alpha_i|\alpha_i)=2$ we must have
$(\alpha_i|\alpha_j)<0$ for some pair $\alpha_i\neq\alpha_j$, hence
$\alpha_i + \alpha_j \in \Delta$ and we can make the combination
shorter, contrary to our assumption. Hence $\Delta = \Phi$. 
Same argument shows that 
if $\Phi$ is of type $B$ or $BC$ then $\Delta$ contains all short
roots, and if $\Phi$ is of type $C$  then $\Delta$ contains all long roots. 

Let $\Phi$ be of type $C$. Then $\Delta$ contains all
long roots of $\Phi$. Since $\Phi$ is a minimal Cartan root system
containing $\Delta$, the latter must contain at least one extra-long
root $\alpha$. Let $\beta$ be a long root such that
$(\alpha|\beta)=-2$. Then by \prop{rank2}(iii) $\Delta$ contains the
whole root system of type $C_2$ generated by $\alpha$ and $\beta$,
therefore the extra-long root $\alpha+2\beta$ also belongs to
$\Delta$. Continuing this argument we get that all extra-long roots lay
in $\Delta$, hence $\Delta=\Phi$. 

Assume now that $\Phi$ is of type $B_n$. When $n=1$ or $2$ we refer
to \prop{rank1} and \prop{rank2}, so assume that $n\ge3$. Let
$\alpha_1,\ldots,\alpha_n\in \Delta\subset \Phi$
be a basis of $\Lambda$ consisting of pairwise orthogonal short
roots. Let $\Phi_l$ be the set of all long roots in $\Phi$. They form
a simply-laced Cartan root system of type $D_n$ ($A_3$ if $n=3$). 
The root system $\Delta$ must contain some long roots too. Let 
$\Delta_l\subset \Phi_l$ be the set of long roots of $\Delta$. They
must form a root system as well. It is not too difficult to see that
for $\Delta$ to be indecomposable  the root system $\Delta_l$ must be
either equal to the whole $\Phi_l$ or to
$\set{\alpha_i-\alpha_j}{i\neq j}$, in which case $\Delta_l$ is of
type $A_{n-1}$.  This choice of $\Delta_l$ is
unique up to the action of the Weyl group. Therefore $\Delta$ is
either equal to $\Phi$ or is of type $B^0_n$. 

Finally, let $\Phi$ be of type $BC_n$. Then by the result of the
previous paragraph the set $\Delta_l$ of long roots must at least
contain the set $\set{\alpha_i-\alpha_j}{i\neq j}$. Also, $\Delta$
must contain at least one extra-long root. So, as in the case of type
$C$, using (iii) of \prop{rank2} we obtain that $\Delta=\Phi$.
\end{proof}

If the root system $\Delta$ is of the type $A$-$D$-$E$, then the
corresponding  conformal algebra
$\goth L = \bigoplus_{\lambda\in\Delta}\goth L_\lambda \subset V_\Lambda$ 
generated by the set $\set{v_\lambda}{\lambda\in\Delta}$ is the affine
Kac-Moody conformal algebra and $V_\Lambda$ is Frenkel-Kac-Segal
construction of its basic representation, see \cite{fk,segal}.
In this case $\goth L$ is a central extension of a simple loop
 algebra, see \sec{affine}.  
If $\Delta$ is of type $C$ then $\goth L$ is also a central extension
of a simple conformal 
algebra, which is a generalization of the algebra 
$\goth K$, constructed in \sec{weyl}.

%% file: finite.tex
\subsection{Finite root systems}
\label{sec:finite}
In this section we describe all possible finite root systems. Let 
$\Delta\subset \Lambda$ be an indecomposable root system and assume
that $|\Delta|<\infty$. As we have seen in \sec{rootsys}, the bilinear
form $(\,\cdot\,|\,\cdot\,)$ on $\Lambda$ must be either positive or 
semi-positive definite. Assume that it is semi-positive definite. 
Let $\pi:\Lambda\to \bar\Lambda$ be the projection of $\Lambda$ onto the
positive definite lattice $\ol\Lambda$, and let 
$\ol\Delta = \pi(\Delta)\ssm \{0\}\subset \ol\Lambda$ be the positive
definite finite root system obtained from 
the projection of $\Delta$, see
\sec{rootsys} for the definitions. The root system $\ol\Delta$
decomposes into a disjoined union 
$\ol\Delta = \ol\Delta_1 \sqcup \ldots \sqcup\ol\Delta_l$ of indecomposable
root systems, each of them must be of one of the types described in
\thm{pd}. Denote $\Delta_i = \pi\inv\big(\ol\Delta_i\big)\cap \Delta$.

If for some $\ol\Delta_i$ we have $\# \pi\inv(\alpha) =1$ for all 
$\alpha \in \ol\Delta_i$ then $\Delta$ decomposes as 
$\Delta_i \sqcup \big(\bigcup_{j\neq i} \Delta_j\big)$, which is a
contradiction. So we assume that each  
$\Delta_i$ is a semi-positive definite root system.

\begin{lemma}\label{lem:finite}\sl
Let $\Delta$ be an indecomposable semi-positive definite root system
such that $\ol\Delta$ is a positive definite indecomposable root
system of type other than $B$ or $B^0$. Then $|\Delta|=\infty$.
\end{lemma}

\begin{proof}
The root system $\Delta$ must contain some isotropic roots, otherwise
it would be positive definite. At least some isotropic root $\delta$
must be of the form $\delta = \alpha+\beta$, where $\alpha$ and
$\beta$ are real roots. If $\Delta$ is not of type $BC$ then $\alpha$
and $\beta$ have square lengths more than 1 and hence we are in the
situation of (v), (vi) or (vii) of \prop{rank2}, so $k\delta \in
\Delta$ for all integer $k$ and $\Delta$ is infinite. If $\Delta$ is
of type $BC$, then it might happen that $\alpha$ and $\beta$ have
length 1. If this is the case, let $\alpha'$ and $\beta'$ be real roots such
that $\ol{\alpha'} = 2\bar \alpha$ and $\ol{\beta'} = 2\bar
\beta$. Then by \prop{rank2}(vii), $\delta' = \alpha' + \beta'$ is an isotropic root such
that $k\delta'$ is also a root for all integer $k$.
\end{proof}

Return now to our finite root system $\Delta$. The lemma implies that
all indecomposable components $\ol\Delta_i$ of $\ol\Delta$ are of type
either $B$ or $B^0$. On the other hand, assume we are given a positive
definite root system $\ol\Delta =  \ol\Delta_1 \sqcup \ldots\sqcup
\ol\Delta_l\subset \ol\Lambda$ such that all components $\ol\Delta_i$ are of type either
$B$ or $B^0$. Let  $\ol\Delta_{\text{s}}$ (respectively, $\ol\Delta_{\text{l}}$) be
the set of short (respectively, long) roots of
$\ol\Delta$. There are many degrees of freedom in
reconstructing the finite semi-positive definite root system
$\Delta$. First we choose an arbitrary lattice $\Lambda_0$ and set 
$\Lambda = \ol\Lambda \oplus \Lambda_0$. Then in each $\ol\Delta_i$ of
type $B$ we choose a subsystem $\ol\Delta_i'$ of type $B_0$. 
 Denote by  $\Omega\subset \ol\Delta_{\text{l}}$ be the set
of all long roots in $\ol\Delta$ which do not get into any of the root systems
$\ol\Delta_i$ or $\ol\Delta'_i$
of type $B^0$. For each $\alpha \in \Omega$ we choose an arbitrary 
isotropic vector $\delta(\alpha) \in \Lambda_0$ such that
$\delta(\alpha) = -\delta(-\alpha)$ and for each short
root $\beta \in \ol\Delta_{\text{s}}$ we choose an arbitrary finite set
$\Sigma(\beta) \subset \Lambda_0$ such that
$\Sigma(\beta)=-\Sigma(\beta)= \Sigma(-\beta)$. We impose the
following restriction: If $\alpha \in \Omega$ and $\beta$ is a short
root such that $(\alpha|\beta)\neq 0$ then $\delta(\alpha) \in
\Sigma(\beta)$. Now we set 
$$
\Delta = \set{\beta+\delta}{\beta \in \ol\Delta_{\text{s}}, \ \ \delta
\in \Sigma(\beta)}
\cup
\set{\alpha+\delta(\alpha)}{\alpha \in \Omega}
\cup
\big(\ol\Delta_{\text{l}}\ssm \Omega\big).
$$

To summarize:

\begin{theorem}\label{thm:finite}\sl
Assume $\Delta\subset \Lambda$ be a finite indecomposable root
system. Then either $\Lambda$ is positive definite and then $\Delta$
is of one of the types described in \thm{pd} or $\Lambda$ is
semi-positive definite, the positive definite quotient of $\Delta$
decomposes into a disjoined union  
$\ol\Delta = \ol\Delta_1 \sqcup \ldots \sqcup\ol\Delta_l$ of finite positive
definite root systems of type either $B$ or $B^0$ and $\Delta$ could
be reconstructed from $\ol\Delta$ by the above procedure.
\end{theorem}

%% file: ears.tex
\subsection{Connection to extended affine root systems}
\label{sec:ears}
In this section we point out the relations with the theory of extended
affine root systems (EARS), see e.g.  \cite{aabgp}.
By \prop{rank2}, for any two vectors $\alpha, \beta \in \Delta$ such that
$(\alpha|\alpha)\neq0$ the Cartan number  
$\<\alpha, \beta\> = \frac{2(\alpha|\beta)}{(\alpha|\alpha)}$ 
is an integer. This already makes $\Delta$ look similar to an EARS. 
To make the similarity even more complete we must impose the
following indecomposability assumption: 
\begin{equation}\label{fl:indec}
\forall\, \delta \in \Delta_0\quad
\exists\, \alpha\in \Delta^\times \quad \text{such that}\quad
\delta + \alpha \in \Delta^\times. 
\end{equation}

It is easy to see that if a root system $\Delta$ satisfies \fl(indec)
and $\ol\Delta$ is indecomposable then $\Delta$ lays inside 
an EARS, as defined in \cite[page 1]{aabgp}.
 The following
theorem shows that in the case when 
there are no short roots, the structure of $\Delta$ is much
simpler than the structure of a general EARS. 

\begin{theorem}\label{thm:spd}\sl
Assume that $\ol\Delta$ is an indecomposable root system of one of the
following types: $A_n$,
$D_n$, $E_n$, $B'_1$ or $C_n$, i.e. $\Delta$ does not contain short roots. 
Assume also that the condition \fl(indec) holds. Then for any 
$\delta \in \Delta_0$ and $\alpha \in \Delta$ we have $\delta + \alpha
\in \Delta$. 
\end{theorem}

The theorem asserts that $\Delta_0$ is a sublattice in $\Lambda_0$ and
for any $\alpha \in \ol\Delta$ the whole equivalence class $\alpha +
\Delta_0$ belongs to $\Delta$. 

\begin{proof}
Let $\goth L = \bigoplus_{\alpha\in \Delta}\goth L_\alpha\subset
V_\Lambda$ be the conformal superalgebra generated by the set
$\set{v_\alpha}{\alpha\in \Delta}$. We claim that for any $h \in \goth
h = \Bbbk \otimes \Lambda$ we have $h(-1)\,v_\delta \in \goth L_\delta$. 

Let us first show that the theorem follows from this claim. 
Let $\delta \in \Delta_0$ and $\alpha \in \Delta$. If $\alpha \in
\Delta^\times$, then using \fl(confjac) and \fl(vanvb), we get
$$
\bigl(\alpha(-1)\,v_\delta\bigr) \ensquare 0
v_\alpha = \pm (\alpha|\alpha)\, v_{\alpha+\delta}\in\goth L,
$$
hence $\alpha +
\delta \in \Delta^\times$. If $\alpha \in \Delta_0$, take some $\beta
\in \Delta^\times$ and then using as before, \fl(confjac) and
\fl(vanvb), we get 
$$
\bigl(\beta(-1)\,v_\delta \bigr)\ensquare 0 \bigl(\beta(-1)\,v_\alpha
\bigr) = \pm (\beta|\beta) \,\delta(-1)\, v_{\alpha+\delta}\in\goth L,
$$ 
hence $\alpha+\delta \in \Delta_0$.

Let us now prove the claim.
The condition \fl(indec) assures that any isotropic root $\delta \in
\Delta_0$ is obtained as a sum $\delta = \alpha + \beta$ of two real
roots $\alpha, \beta \in \Delta^\times$. Since $\Delta$ does not have
any short roots, the pair $\alpha, \beta$ must generate a root system
of type either (v) (vi) or (vii) of \prop{rank2}. So \prop{rank2}
implies that $\alpha(-1)\,v_\delta, \beta(-1)\,v_\delta\goth \in L_\delta$.

Assume now that $\lambda \in \Delta^\times$ is such that 
$\lambda(-1)\,v_\delta \in \goth L_\delta$ and $\mu \in \Delta^\times$ 
satisfies $(\lambda|\mu)\neq 0$. Then we have
$$
\bigl(\lambda(-1)\,v_\delta\bigr) \ensquare 0 v_\mu = \pm(\lambda|\mu)\,
v_{\mu+\lambda}\in \goth L, 
$$ 
hence $\mu+\delta \in \Delta^\times$. Therefore the real roots $\mu$
and $\mu+\delta$ form a root system of type (v), (vi) or (vii) of
\prop{rank2}, so we get that $\mu(-1)\, v_\delta \in \goth L_\delta$. 

It follows that for every real root $\lambda\in \Delta^\times$ which is
not orthogonal to either $\alpha$ or $\beta$ we have
$\lambda(-1)\,v_\delta \in\goth L_\delta$, therefore, since
$\ol\Delta$ is indecomposable, $\lambda(-1)\,v_\delta \in\goth
L_\delta$ for all $\lambda \in \Delta^\times$.  The condition
\fl(indec) implies that $\goth h = \op{Span}_\Bbbk \Delta^\times$, and
the claim follows.
\end{proof}